\documentclass[preprint, 11pt]{elsarticle}
\usepackage{lineno}
\pdfoutput=1

\usepackage{amsmath}
\usepackage{amssymb}
\usepackage{mathtools}
\usepackage{mathrsfs}
\usepackage{color}
\usepackage{xcolor}
\usepackage{graphicx}
\usepackage[percent]{overpic}
\usepackage{url}
\usepackage{caption}
\usepackage{subcaption}
\usepackage{enumerate}
\usepackage{overpic}
\usepackage{amsthm}
\usepackage{moreverb}
\usepackage{multirow}
\usepackage{multicol}
\usepackage[pdftex,colorlinks,bookmarksopen,bookmarksnumbered,citecolor=red,urlcolor=red]{hyperref}
\usepackage{mwe}
\usepackage{graphbox}
\usepackage{textcomp}
\usepackage{tabularx}
\usepackage{algorithm}
\usepackage{algorithmicx}
\usepackage[noend]{algpseudocode}

\def\n{{\bm n}}

\def\0{\boldsymbol{0}}

\def\xbf{\boldsymbol{x}}

\def\calcf{\mathbf{\mathcal{C}}}
\def\calsf{\mathbf{\mathcal{S}}}
\def\caldf{\mathbf{\mathcal{D}}}
\def\caluf{\mathbf{\mathcal{U}}}
\def\calsigf{\mathbf{\bm{\Sigma}}}
\def\calvf{\mathbf{\mathcal{V}}}

\def\calmf{\mathbf{\mathcal{M}}}
\def\calpf{\mathbf{\mathcal{P}}}

\def\caltf{\mathbf{\mathcal{T}}}

\def\tildeq{\widetilde{q}}
\def\tildep{\widetilde{\psi}}

\def\cl {\nonumber \\}

\newcommand{\bm}[1]{\mbox{\boldmath{$#1$}}}
\DeclarePairedDelimiter\abs{\lvert}{\rvert}%
\definecolor{ForestGreen}{RGB}{34,139,34}
\newcommand\norm[1]{\lVert#1\rVert}
\def\div{\nabla\cdot}
\def\grad{\nabla}

\newcommand{\anna}[2][cyan]{{\textcolor{#1}{#2}}}

\usepackage{verbatim}
\usepackage{graphicx}
\usepackage{epstopdf}

\newtheorem*{remark}{Remark}

\begin{document}

\begin{frontmatter}


\title{Randomized Proper Orthogonal Decomposition for data-driven reduced order modeling of a two-layer quasi-geostrophic ocean model}





\author[Houston]{Lander Besabe}
\ead{lybesabe@central.uh.edu}

\author[UNIPA]{Michele Girfoglio}
\ead{michele.girfoglio@unipa.it}

\author[Houston]{Annalisa Quaini\corref{mycorrespondingauthor}}
\ead{aquaini@central.uh.edu}
\cortext[mycorrespondingauthor]{Corresponding author}

\author[SISSA]{Gianluigi Rozza}
\ead{grozza@sissa.it}

\address[Houston]{Department of Mathematics, University of Houston, 3551 Cullen Blvd, Houston TX 77204, USA}

\address[UNIPA]{University of Palermo, Department of Engineering, Viale delle Scienze, Ed. 7, Palermo, 90128, Italy}

\address[SISSA]{SISSA, International School for Advanced Studies, Mathematics Area, mathLab, via Bonomea, Trieste 265 34136, Italy}

\begin{abstract}
The two-layer quasi-geostrophic equations (2QGE) serve as a simplified model for simulating wind-driven, stratified ocean flows. However, their numerical simulation remains computationally expensive due to the need for high-resolution meshes to capture a wide range of turbulent scales.
This becomes especially problematic when several simulations need to be run because of, e.g., uncertainty in the parameter settings. To address this challenge, 
we propose a data-driven reduced order model (ROM) for the 2QGE that leverages randomized proper orthogonal decomposition (rPOD) and long short-term memory (LSTM) networks. To efficiently generate the snapshot data required for model construction, we apply a nonlinear filtering stabilization technique
that allows for the use of
larger mesh sizes compared to a direct numerical simulations (DNS). Thanks to the use of rPOD to extract the dominant modes from the snapshot matrices, we achieve up to 700 times speedup over the use of deterministic POD.
LSTM networks are trained with the modal coefficients associated with the snapshots to enable the prediction of the time- and parameter-dependent modal coefficients during the online phase, which is hundreds of thousands of time faster than a DNS. 
We assess the accuracy and efficiency of our rPOD-LSTM ROM through an extension of
a well-known benchmark called double-gyre wind forcing test. The dimension of the parameter space in this test is increased from two to four. 
\end{abstract}


\end{frontmatter}

\section{Introduction} \label{sec:intro}
The quasi-geostrophic equations (QGE) are a widely used model for geophysical flows, capturing critical dynamics such as baroclinic instabilities and large-scale vortices for ocean and atmospheric flows. Here, we focus on their use for ocean modeling.
The QGE model, considered simpler for numerical and analytical investigation \cite{Hoskins1975}, 
represents the ocean as one layer with uniform depth, density, and temperature. 
The  two-layer quasi-geostrophic equations (2QGE) are an attempt to capture the complexity of stratification 
by adding a second dynamically active layer. 
Although still regarded as a simplified model, the simulation of 
ocean dynamics through the 2QGE poses significant challenges. 
This paper aims at addressing two of such challenges: the first is related to the fact that the 2QGE model features several parameters, whose exact value might be only known to vary in given ranges, 
and the second is the high computational cost for each simulation once the parameters are set. 


ROMs (see, e.g., \cite{peter2021modelvol1,benner2020modelvol2,benner2020modelvol3,hesthaven2016certified, malik2017reduced,rozza2008reduced} for reviews) have emerged as the methodology of choice to reduce the computational cost when traditional flow simulations (Full Order Models, FOMs) have to be carried out for several parameter values, as in the case of uncertain parameter values. 
ROMs replace the FOM with a lower-dimensional approximation that captures the essential flow behavior. The reduction in computational time is achieved through a two-step procedure. In the first step, called {offline}, one constructs a database of several FOM solutions associated to given times and/or physical parameter values. The FOM database is used to generate a reduced basis, which is (hopefully much) smaller than the high-dimensional FOM basis but still preserves the essential flow features. 
Techniques such as proper orthogonal decomposition (POD) \cite{Berkooz1993} and dynamic mode decomposition (DMD) \cite{Schmid2010} are commonly used to perform this step. 
In the second step, called online, 
one uses this reduced basis to quickly compute the solution for newly specified times and/or parameter values.
ROMs have become popular across different scientific disciplines, such as structural mechanics \cite{Casenave2019, Guo2018, Lin1995, Oliver2017} 
and biomedical engineering \cite{Lucas2019, Pfaller2020, Pfaller2022, Itu2012, Zhang2020}.

Although performed once, the offline stage of a ROM could have a significant computational cost for two reasons: i) one needs to compute several FOM solutions to construct a database and ii) one needs to generate a reduced basis from a possibly
large database. The cost of each FOM solution for the 2QGE can be especially high
because impractically fine meshes are needed to resolve the full spectra of turbulence. To make the computational cost manageable, researchers typically resort to lower resolution
meshes and the use of the so-called eddy viscosity parameterization, i.e., the viscosity of water (order of $10^{-6}$ m$^{2}$/s) is replaced with an artificially large viscosity coefficient
(order of few thousand \cite{Bryan1963,Gates1968}
to a few hundred \cite{Holland1975,Berloff1999,BERLOFF_KAMENKOVICH_PEDLOSKY_2009,Tanaka2010})
to compensate for the diffusion
mechanisms that are not captured due to mesh under-resolution. In this paper, 
we resort to a nonlinear filtering
stabilization technique we recently introduced \cite{Besabe2024} to further
reduce the computational cost for a given coarse mesh, so that each FOM solution becomes cheaper.

In the recent years, researchers have looked into randomized linear algebra methods to alleviate the cost 
associated with the above point ii). 
A randomized DMD algorithm based on 
randomized singular value decomposition (SVD) was introduced in \cite{Erichson2019}
and assessed with two test cases, i.e, 2D flow behind a cylinder and 
world climate data, with time as the 
only parameter. It was shown that 
randomized DMD is 
up to 60 times faster than deterministic DMD on a CPU, while preserving good accuracy. 
An adaptive randomized DMD algorithm was explored in \cite{Bistrian2018} for the development of ROMs for the shallow water equations. It was found that this  
adaptive randomized DMD algorithm is substantially faster than both classic and energetic/optimized DMD \cite{Tissot2014}, i.e., DMD with modes that are chosen based on their amplitude, frequency, or growth rate. In \cite{Buhr2018}, adaptive randomized rangefinder algorithms were used to efficiently extract the dominant modes from 
local solutions of the Helmholtz equation and a linear elasticity problem. 
In \cite{Saibaba2020}, 
randomized versions of discrete empirical interpolation method (DEIM) were developed, analyzed, and tested on numerical examples from \cite{Peherstorfer2014}, which are nonlinear interpolation problems where the reference functions exhibit sharp peaks, 
and advection-diffusion problems with three varying parameters. 
In \cite{Yu2020}, balanced POD \cite{Rowley2005} was performed efficiently through randomized methods for both the primal and adjoint systems of large-scale linear systems.
The numerical tests considered in \cite{Yu2020} are heat transfer and advection-diffusion problems, with time as the only parameter. 
Randomized POD, DMD, and DEIM based on randomized SVD were numerically tested 
in \cite{Alla2019}. They have shown to be more efficiently than their deterministic counterparts for 
linear and semi-linear problems, with at most two varying parameters. 
ROMs using randomized POD were developed in \cite{Rajaram2020} and tested for a nonlinear Poisson problem 
and the transonic RAE 2822 airfoil problem. 
In \cite{Bach2019}, different low-rank approximation methods (e.g., incremental and randomized SVD) were considered to develop ROMs for large nonlinear problems. These methods were tested with different problems: the 1D Burgers' problem, the Taylor bar impact test, modeling high-velocity impact of a circular bar onto a rigid surface, 
and a problem in the automotive industry involving contact, buckling, plasticity, and large deformations. 

In this work, we develop a data-driven ROM (i.e., the online stage is blind to the underlying model and relies only on solution data from the FOM) for the 2QGE that uses randomized POD to generate the reduced basis. 
To find the coefficients for the reduced basis functions
and thus obtain
the ROM approximation of each variable, we adopt 
a simple LSTM architecture \cite{Besabe2024b,Rahman-2019}.
We test the accuracy and efficiency of the proposed rPOD-LSTM ROM
as the dimensionality of the parameter space is increased up to four (time plus three physical coefficients). 
We show that our ROM is computationally efficient, even in the case of high-dimensional parameter spaces, while maintaining predictive accuracy for the time-averaged variables of interest, with relative $L^2$ errors ranging from 1E-02 to 1E-01 when the parameter dimension is up to three. 
In the case of a four-dimensional parameter space, the relative $L^2$ errors increase. However, they can be contained with a combination of finer parameter sampling, which can be easily handled thanks to the randomized POD and ROM closure (see, e.g., \cite{Girfoglio2023,fluids9080178}).

The rest of the paper is organized as follows. Sec.~\ref{sec:2qge} describes the 2QGE model and Sec. \ref{sec:fom} presents the FOM. Sec.~\ref{sec:rom} discusses the building blocks of our ROM: randomized POD  and LSTM networks. Numerical results are presented in
Sec.~\ref{sec:num_res}. Finally, 
conclusions are drawn in Sec.~\ref{sec:conclusions}.

\section{Governing Equations}\label{sec:2qge}

We consider an ocean with two layers with uniform depth, density, and temperature, 
on a two-dimensional rectangular domain $\Omega = [x_0, x_f]\times [-L/2, L/2]$. 
We will call layer 1 the top layer, i.e., the layer driven by the wind, and 
layer 2 the bottom layer. We assume 
that the meridional length $L$ is substantially larger than the depths of both layers, denoted by $H_1$ and $H_2$. We refer the reader to, e.g., \cite{Marshall1997,Chassignet1998,Berloff1999,DiBattista2001,BERLOFF_KAMENKOVICH_PEDLOSKY_2009}, for a thorough description of other assumptions leading to the derivation of 2QGE from the Navier-Stokes equations.

Let $(0,T]$ be the time interval of interest. The non-dimensional formulation of 2QGE reads: find potential vorticities $q_l$ and stream functions $\psi_l$, for $l=1,2$, such that 
\begin{align}
    &\frac{\partial q_1}{\partial t} - \nabla \cdot \left(\left(\nabla \times \bm{\Psi}_1 \right)q_1\right) + \frac{Fr}{Re~\delta}\Delta(\psi_2-\psi_1) - \frac{1}{Re}\Delta q_1 = F, \label{eq:qge2_1}\\
    &\frac{\partial q_2}{\partial t} - \nabla \cdot \left(\left(\nabla \times \bm{\Psi}_2 \right)q_2\right) + \frac{Fr}{Re~(1-\delta)}\Delta(\psi_1-\psi_2) + \sigma\Delta\psi_2- \frac{1}{Re}\Delta q_2 = 0, \label{eq:qge2_2}\\
    &q_1 = Ro\Delta\psi_1 + y + \frac{Fr}{\delta}\left(\psi_2 - \psi_1 \right), \label{eq:kin_eq1} \\
    &q_2 = Ro\Delta\psi_2 + y + \frac{Fr}{1-\delta}\left(\psi_1 - \psi_2 \right), \label{eq:kin_eq2}
\end{align}
in $\Omega\times(0,T)$, where $\bm{\Psi}_l=(0, 0, \psi_l)$, $F$ is the wind forcing on the top layer, $y$ is the non-dimensional vertical coordinate of the domain, $\delta = \frac{H_1}{H_1+H_2}$ denotes the aspect ratio of the layer depths, and $\sigma$ is the friction coefficient at the bottom of the ocean. Problem \eqref{eq:qge2_1}-\eqref{eq:kin_eq2} features also the following non-dimensional numbers: 
the Froude number $Fr$, 
the Reynolds number $Re$, and
the Rossby number $Ro$. To give
their definitions, let $U$ be the characteristic velocity scale, 
$H=H_1 + H_2$ the total ocean depth, and
$g' = g\Delta\rho/\rho_1$ the reduced gravity, with $g$ the gravitational constant, $\Delta\rho$ the density difference between the two layers, and $\rho_1$ is the density of the top layer.
Moreover, 
let $\nu$ be the (constant) eddy viscosity coefficient and $\beta$ be the gradient of the Coriolis frequency $f$, i.e.,  $f \approx f_0 + \beta y$, where $f_0$ is the local Earth rotation rate at $y = 0$. 
Then:
\begin{equation*}
    Fr = \frac{f_0^2U}{g'\beta H}, \quad Re=\frac{UL}{\nu},\quad Ro = \frac{U}{\beta L^2}.
\end{equation*}
The reader interested in details about the the non-dimensionalization of the 2QGE is referred to, e.g., \cite{San2012, Salmon1978, Medjo2000, Fandry1984, Mu1994}.

\begin{remark}
 We note that we call $\psi_i$, $i = 1, 2$, stream function although it is actually stream function with the opposite sign. Consequently, velocity is $-\nabla \times \bm{\Psi}_i$. Although this is a potential source of confusion, 
 it is the terminology used in, e.g., \cite{San2012, Fandry1984, Ikeda1981, Zalesny2022} 
 and thus we stick to it.
\end{remark}

To complete the model, we follow \cite{Cummins1992, Ozgok1998, San2012} by prescribing free-slip and impenetrable boundary conditions and start the system from a rest state, i.e., we impose the following boundary and initial conditions:
\begin{align} \label{qge-bdry}
    \psi_l = 0, &\quad \mbox{on }\partial\Omega\times(0,T),\\
    q_l = y, &\quad \mbox{on }\partial\Omega\times(0,T),\\
    q_l(x,y) = q_l^0 = y, &\quad \mbox{in }\Omega\times\{0\},
\end{align}
for $l=1,2$.

\section{Full Order Method} \label{sec:fom}

With a slight abuse of notation, we call 
direct numerical simulation (DNS) a simulation with mesh size $h$ that is smaller than the Munk scale defined by
\begin{equation} \label{eq:Munk}
    \delta_M = L\sqrt[3]{\frac{Ro}{Re}},
\end{equation}
where the abuse comes from the fact that
$Re$ is defined with the eddy viscosity coefficient, instead of the actual viscosity of water. Even when using an artificially
large viscosity coefficient, the requirement $h < \delta_M$ is typically very restrictive
and leads to simulations that can take days to complete. See, e.g., 
\cite{Besabe2024}. Thus, we prefer to adopt an alternative FOM to reduce the time needed
to collect the snapshots.

To enable the use of coarser (and thus computationally cheaper) meshes with $h>\delta_M$, we add a nonlinear differential low-pass filter to the 2QGE \cite{Besabe2024}. Specifically, we consider
the Helmholtz filter given by: 
\begin{equation} \label{eq:filter}
    -\alpha^2\div\left(a(q_l)\grad \overline{q}_l\right) + \overline{q}_l = q_l,
\end{equation}
where $\alpha>0$ is the \textit{filtering radius} 
and $0<a(\cdot)\leq1$ is a scalar function called an \textit{indicator}, which satisfies the following conditions:
\begin{align*}
    a(\cdot)\simeq 0 & \mbox{ in regions where the flow field requires little-to-no regularization;}\\
    a(\cdot)\simeq 1 & \mbox{ in regions where the flow field requires } O(\alpha) \mbox{ regularization.}
\end{align*}
In our previous work \cite{Besabe2024},
we showed that a suitable indicator function is
\begin{equation} \label{eq:nl-ind_fun}
    a(q) = \frac{\vert\grad q\vert}{\norm{\grad q}_\infty}.
\end{equation}

Adding the filter \eqref{eq:filter}-\eqref{eq:nl-ind_fun} to the 2QGE \eqref{eq:qge2_1}-\eqref{eq:kin_eq2}, we
obtain the following problem: 
find potential vorticities $q_l$, filtered potential vorticities $\overline{q_l}$, and stream functions $\psi_l$, for $l=1,2$, such that 
\begin{align}
    &\frac{\partial q_1}{\partial t} - \nabla \cdot \left(\left(\nabla \times \bm{\Psi}_1 \right)q_1\right) + \frac{Fr}{Re~\delta}\Delta(\psi_2-\psi_1) - \frac{1}{Re}\Delta q_1 = F, \label{eq:fil1}\\
    &-\alpha^2_1\div\left(a(q_1)\grad \overline{q}_1\right) + \overline{q}_1 = q_1\label{eq:fil2}\\
    &Ro\Delta\psi_1 + y + \frac{Fr}{\delta}\left(\psi_2 - \psi_1 \right) = \overline{q}_1, \label{eq:fil3} \\
    &\frac{\partial q_2}{\partial t} - \nabla \cdot \left(\left(\nabla \times \bm{\Psi}_2 \right)q_2\right) + \frac{Fr}{Re~(1-\delta)}\Delta(\psi_1-\psi_2) + \sigma\Delta\psi_2- \frac{1}{Re}\Delta q_2 = 0, \label{eq:fil4}\\
    &-\alpha^2_2\div\left(a(q_2)\grad \overline{q}_2\right) + \overline{q}_2 = q_2 \label{eq:fil5}\\
    &Ro\Delta\psi_2 + y + \frac{Fr}{1-\delta}\left(\psi_1 - \psi_2 \right) = \overline{q}_2, \label{eq:fil6}
\end{align}
in $\Omega\times(0, T)$. We call this system the 2QG-NL-$\alpha$ model. 
The artificial diffusion introduced by 
the filter is proportional to $\alpha_l$
in layer $l$. For this reason, we allow to set
different values of the filtering radius in the two
layers. For instance, one may choose $\alpha_2 < \alpha_1$
since the equations governing the dynamics in the bottom layer feature an extra diffusive term coming from 
friction with the bottom of the ocean.
For a more detailed discussion on this method, see \cite{Besabe2024, Girfoglio_JCAM2023}. 

To further contain the computational cost at the FOM level, we use a segregated algorithm. Let us denote by $f^n$ the approximation of a quantity $f$ at time $t_n = n\Delta t$, where $\Delta t = T/N$ for some $N\in\mathbb{N}$.
At time $t^{n+1}$, the algorithm reads: given $(q_l^n,\overline{q}_l^n, \psi_l^n)$, for $l=1,2$, perform
\begin{itemize}
    \item[-] Step 1: find the potential vorticity of the top layer $q_1^{n+1}$ such that
    \begin{align}\label{eq:seg-alg1}
        \frac{1}{\Delta t}q_1^{n+1} - \div\left(\left(\nabla \times \bm{\Psi}_1^{n}\right)q_1^{n+1}\right) - \frac{1}{Re}\Delta q_1^{n+1} = b_1^{n+1} 
        - \frac{Fr}{Re~\delta} \Delta \left( \psi_2^{n} - \psi_1^{n} \right).
    \end{align}
    \item[-] Step 2: find the filtered potential vorticity for the top layer $\bar{q}_1^{n+1}$  such that
    \begin{equation} \label{eq:seg-alg2}
        -\alpha^2_1 \nabla \cdot \left( a_1^{n+1}\nabla \bar{q}_1^{n+1}\right) + \bar{q}_1^{n+1} = q_1^{n+1}.
    \end{equation}
    \item[-] Step 3: find the stream function of the top layer $\psi_1^{n+1}$ such that:
    \begin{equation} \label{sf-1-discreet}
         Ro\Delta \psi_1^{n+1} + y - \frac{Fr}{\delta} \psi_1^{n+1} = \bar{q}_1^{n+1} - \frac{Fr}{\delta}\psi_2^{n}.
    \end{equation}
    \item[-] Step 4: find the potential vorticity  of the bottom layer $q_2^{n+1}$ such that:
    \begin{align}
        &\frac{1}{\Delta t}q_2^{n+1} - \div\left(\left(\nabla \times \bm{\Psi}_2^{n}\right)q_2^{n+1}\right) - \frac{1}{Re}\Delta q_2^{n+1}= b_2^{n+1} -\sigma \Delta \psi_2^{n} + \frac{Fr}{Re~(1 - \delta)} \Delta \left( \psi_1^{n+1} - \psi_2^{n} \right).
    \end{align}   
    \item[-] Step 5: find the filtered potential vorticity for the bottom layer $\bar{q}_2^{n+1}$ such that:
    \begin{equation} \label{eq:seg-alg5}
        -\alpha^2_2 \nabla \cdot \left( a_2^{n+1}\nabla \bar{q}_2^{n+1}\right) + \bar{q}_2^{n+1} = q_2^{n+1}.
    \end{equation}
    \item[-] Step 6: find the stream function of the bottom layer $\psi_2^{n+1}$ such that:
    \begin{equation}\label{eq:seg-alg6}
         Ro\Delta \psi_2^{n+1} + y - \frac{Fr}{1-\delta} \psi_2^{n+1} = \bar{q}_2^{n+1} - \frac{Fr}{1-\delta}\psi_1^{n+1} .
    \end{equation}
\end{itemize}
We note that we have used Backward Difference Formula of order 1 for the approximation of the time derivatives in 
\eqref{eq:fil1}-\eqref{eq:fil6}.

Regarding the space discretization of \eqref{eq:seg-alg1}-\eqref{eq:seg-alg6}, we employ a Finite Volume (FV) approximation.
For this purpose, we partition domain $\Omega$ into $N_C$ control volumes denoted by $\Omega_k$, $k=1,\dots,N_C$, such that the $\Omega_k$ are pairwise disjoint. Let $\textbf{A}_j$ represent the surface vector associated with each face of the control volume $\Omega_k$, i.e., $\textbf{A}_j = A_j \n_j$ where $A_j$ is the surface area of face $j$ and $\n_j$ its outward unit normal. The FV discretization is derived by integrating the semi-discrete 2QG-NL-$\alpha$ 
model over the computational domain, using the Gauss divergence theorem, and discretizing the integrals. 
To write the resulting set of equations, we denote with 
$q_{l,k}$ and $\psi_{l,k}$ are the average potential vorticity and average stream function  of layer $l$ over the control volume $\Omega_k$.  
Let $q_{l,k}^{j}$ be the potential vorticity associated to the centroid of the $j$-th face and normalized by the volume of $\Omega_k$, which we compute through a linear interpolation scheme over neighboring cells which is second-order accurate.
We approximate the diffusive terms with a second-order centered difference scheme, as well.
For the convective terms, we define
\begin{equation*}
    \varphi_{l,j}=\left(\nabla \times \bm{\Psi}_{l,j}\right)\cdot\textbf{A}_j \text{ with } \bm{\Psi}_{l,j} = (0,0,\psi_{l,j})^T.
\end{equation*}
Then, the scheme discretized in space and time can be written as:
\begin{align}
    &\frac{1}{\Delta t}q_{1,k}^{n+1} - \sum_j\varphi_{1,j}^{n}q_{1,k}^{n+1,j} - \frac{1}{Re}\sum_j \left(\nabla q_{1,k}^{n+1}\right)_j\cdot \textbf{A}_j = b_{1,k}^{n+1} - \frac{Fr}{Re~\delta}\sum_j \left(\nabla(\psi_{2,k}^n - \psi_{1,k}^n)\right)_j\cdot \textbf{A}_j \\
    & -\alpha^2_1 \sum_j a_{1,k}^{n+1} \left(\nabla \bar{q}_{1,k}^{n+1}\right)_j\cdot \textbf{A}_j + \bar{q}_{1,k}^{n+1}=q_{1,k}^{n+1},\\
    & Ro\sum_j\left(\nabla \psi_{1,k}^{n+1}\right)_j\cdot\textbf{A}_j + y_k + \frac{Fr}{\delta}\left(\psi_{2,k}^n-\psi_{1,k}^{n+1}\right) = \bar{q}_{1,k}^{n+1},\\
    & \frac{1}{\Delta t}q_{2,k}^{n+1} - \sum_j\varphi_{2,j}^{n}q_{2,k}^{n+1,j} - \frac{1}{Re}\sum_j \left(\nabla q_{2,k}^{n+1}\right)_j\cdot \textbf{A}_j = b_{2,k}^{n+1} \cl 
    & \qquad \quad + \left(\frac{Fr}{Re~(1-\delta)} - \sigma\right) \sum_j (\nabla\psi_{2,k}^n)_j\cdot \textbf{A}_j - \frac{Fr}{Re~(1-\delta)}\sum_j \left(\nabla\psi_{1,k}^{n+1}\right)_j\cdot \textbf{A}_j,
    \end{align}
    \begin{align}
    & -\alpha^2_2 \sum_j a_{1,k}^{n+1} \left(\nabla \bar{q}_{1,k}^{n+1}\right)_j\cdot \textbf{A}_j + \bar{q}_{1,k}^{n+1}=q_{1,k}^{n+1},\\
     & Ro\sum_j\left(\nabla \psi_{2,k}^{n+1}\right)_j\cdot\textbf{A}_j + y_k - \frac{Fr}{1-\delta}\psi_{2,k}^{n+1}= \bar{q}_{2,k}^{n+1} - \frac{Fr}{1-\delta}\psi_{1,k}^{n+1}, \label{eq:disc-psi2}
\end{align}
for each control volume $\Omega_k$, where 
$b_{l,k}^{n+1}$ is the average discrete forcing and $y_k$ is the vertical coordinate of the centroid of $\Omega_k$.

It has been shown in \cite{Besabe2024} that the
2QG-NL-$\alpha$ model and the corresponding numerical scheme yield an accurate solution for an appropriate choice of $\alpha$, i.e., $\alpha = O(h)$, while greatly reducing the computational cost.

All the FOM simulations are executed using GEA \cite{GEA, GirfoglioFVCA10}, an open-source software package which is built upon the C++ finite volume library OpenFOAM\textsuperscript{\textregistered} \cite{Weller1998}.

\section{Reduced Order Modeling} \label{sec:rom}


Let $\bm{\mu}$ be a vector containing  all the physical parameters of interest belonging to parameter space $\caldf\subset \mathbb{R}^d$. We assume that the solution of the system \eqref{eq:fil1}-\eqref{eq:fil6} can be approximated as follows:
\begin{align}
    q_l(t,\xbf, \bm{\mu})\approx q_{l}^r(t,\xbf, \bm{\mu}) &= \tildeq_l^0(\xbf) + \sum_{i=1}^{N_{q_l}^r} \alpha_{l,i}(t, \bm{\mu})\varphi_{l,i}(\xbf),\label{eq:q_approx} \\
    \psi_l(t,\xbf, \bm{\mu})\approx \psi_{l}^r(t, \xbf, \bm{\mu}) &= \tildep_l^0(\xbf) + \sum_{i=1}^{N_{\psi_l}^r} \beta_{l,i}(t, \bm{\mu})\xi_{l,i}(\xbf), \label{eq:psi_approx}
\end{align}
where $q_{l}^r$ and $\psi_{l}^r$ are the reduced order approximations and $\tildeq_l^0$ and $\tildep_l^0$ are zero-th order approximations of the time-averaged potential vorticity $\tildeq_l$ and stream function $\tildep_l$:
\begin{equation}\label{eq:time_av}
    \tildeq_l(\xbf, \bm{\mu}) = \frac{1}{N^t}\sum_{p=1}^{N^t} q_l(t_p,\xbf, \bm{\mu}), \quad \tildep_l (\xbf, \bm{\mu}) = \frac{1}{N^t}\sum_{p=1}^{N^t} \psi_l(t_p,\xbf, \bm{\mu}).
\end{equation}
See Sec.~\ref{sec:pod} for details on how $\tildeq_l^0$ and $\tildep_l^0$ are computed. 

Notice that in \eqref{eq:q_approx}-\eqref{eq:psi_approx} we are approximating the fluctuations of vorticity and stream function as a linear combination of a few global spatial basis functions with parameter-dependent coefficients.
Sec.~\ref{sec:pod} discusses a process frequently used to 
compute the basis functions $\varphi_{l,i}$ in \eqref{eq:q_approx} and $\xi_{l,i}$ in \eqref{eq:psi_approx} 
and determine the cardinalities $N_{q_l}^r$ and $N_{\psi_l}^r$
of the reduced bases. Sec.~\ref{sec:rpod} presents a technique that accomplishes the same goal more efficiently. Lastly, Sec.~\ref{sec:lstm} covers the procedure to compute for the modal coefficients $\alpha_{l,i}$ and $\beta_{l,i}$ in \eqref{eq:q_approx}-\eqref{eq:psi_approx}.

\subsection{Proper Orthogonal Decomposition} \label{sec:pod}

To compute the global spatial basis functions $\varphi_{l,i}$ and $\xi_{l,i}$ of the reduced spaces in \eqref{eq:q_approx}-\eqref{eq:psi_approx}, proper orthogonal decomposition (POD) is frequently used. Let us briefly summarize the procedure that yields the reduced basis. 

We consider $N^t$ sample time instants $t_p$, $p=1,\dots,N^t$, and $M$ sample parameter vectors $\bm{\mu}_k$, $k=1,\dots,M$. We compute the fluctuations from time-averaged fields by
\begin{align} 
    q_l' (t_p, \xbf, \bm{\mu}_k) &= q_l(t_p,\xbf, \bm{\mu}_k) - \tildeq_l(\xbf, \bm{\mu}_k), \label{eq:fluct_q} \\
    \psi_l' (t_p, \xbf, \bm{\mu}_k) &= \psi_l(t_p, \xbf, \bm{\mu}_k) - \tildep_l(\xbf, \bm{\mu}_k), \label{eq:fluct_psi}
\end{align}
where the time-averaged fields are computed by plugging
$\bm{\mu}_k$ into $\bm{\mu}$ in \eqref{eq:time_av}.
We collect the fluctuation snapshots \eqref{eq:fluct_q}-\eqref{eq:fluct_psi} for all parameters $\bm{\mu}_k$ and time instants $t_p$ and store them in the matrix $\calsf_\Phi\in\mathbb{R}^{N_C\times N^s}$, with $N^s = M\cdot N^t$:
\begin{equation}
    \calsf_\Phi = \left[\Phi(t_1,\bm{x},\bm{\mu}_1)\hspace{0.25em}\cdots \hspace{0.25em} \Phi(t_{N^t},\bm{x},\bm{\mu}_1)\hspace{0.5em} \Phi(t_1,\bm{x},\bm{\mu}_2)\hspace{0.25em}\cdots\hspace{0.25em} \Phi(t_{N^t},\bm{x},\bm{\mu}_M)\right]
\end{equation}
where $\Phi\in\{q_l',\psi_l'\}, l=1,2$. Additionally, once we have sampled the parameter space, we can compute 
the zero-th order approximations in \eqref{eq:q_approx}-\eqref{eq:psi_approx}:
\begin{equation}\label{eq:time_av_zero}
    \tildeq_l^0(\xbf) = \frac{1}{N^t}\sum_{p=1}^{N^t} q_l(t_p,\xbf, \bm{\mu}_c), \quad \tildep_l^0 (\xbf) = \frac{1}{N^t}\sum_{p=1}^{N^t} \psi_l(t_p,\xbf, \bm{\mu}_c),
\end{equation}
where $\bm{\mu}_c$ is the sample in the $\{\bm{\mu}_k\}_{k=1}^{M}\subset\caldf$ closest to a parameter vector of interest $\bm{\mu}$ in terms of the Euclidean distance.

Let $K=\min\{N_C, N^s\}$. Singular value decomposition (SVD) is applied to $S_\Phi$ to obtain
\begin{equation} \label{eq:svd}
    \calsf_\Phi = \caluf_\Phi\calsigf_\Phi\calvf_\Phi^T
\end{equation}
where $\caluf_\Phi \in \mathbb{R}^{N_C\times K}$ and $\calvf_\Phi \in \mathbb{R}^{N^s \times K}$ whose columns are the left and right singular vectors of the snapshot matrix $\calsf_\Phi$, and $\calsigf_\Phi \in \mathbb{R}^{K \times K}$ is the diagonal matrix containing the singular values $\sigma_\Phi^i$ of $\calsf_\Phi$ arranged in descending order. 
The POD basis functions are the first $N_\Phi^r$ left singular vectors in $\caluf_\Phi$, where $N_\Phi^r$ is chosen such that for a
user-provided threshold $0<\delta_\Phi < 1$:
\begin{equation} \label{eq:energy}
    \frac{\sum_{i=1}^{N_\Phi^r} \sigma_\Phi^i}{\sum_{i=1}^{K} \sigma_\Phi^i} \geq \delta_\Phi.
\end{equation}
Eq.~\eqref{eq:energy} means that by selecting the first $N_\Phi^r$ POD modes as basis functions
we retain a user-defined fraction ($\delta_\Phi$) 
of the singular value energy of the system, where we use the word ``energy" loosely. 
The POD basis is considered efficient if $N_\Phi^r\ll K$, that is, we can adequately capture the dynamics of the system with a small number of POD modes.

\begin{figure}
    \centering
    \includegraphics[width=0.5\linewidth]{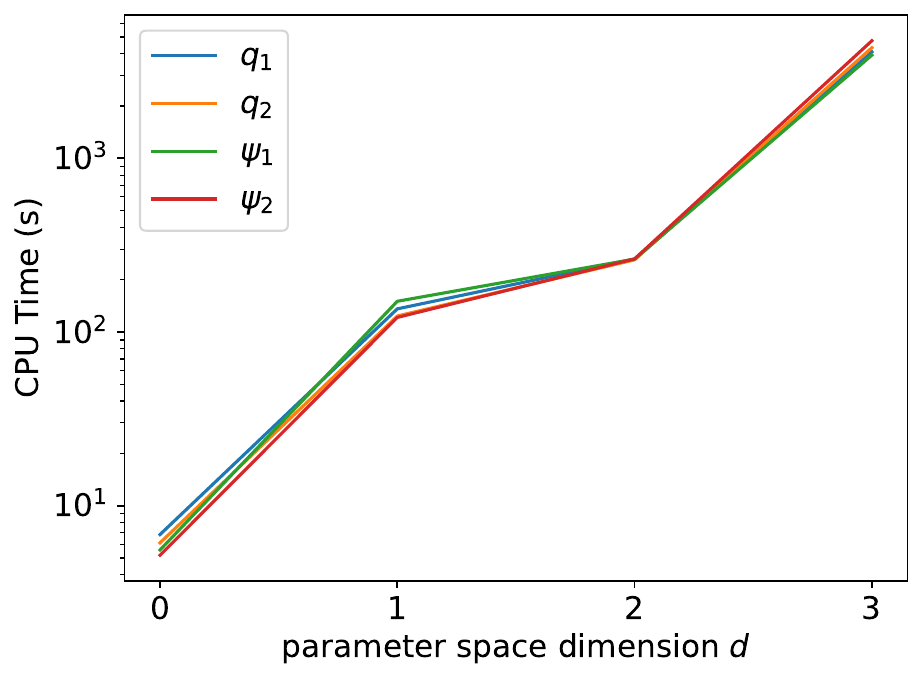}
    \caption{Computational time required to apply POD to a snapshot matrix associated to each variable  for dimension of the physical parameter space $d$ ranging from 0 (i.e., time is the only parameter) to
    3 (i.e., $\delta$, $\sigma$, and $Fr$, in addition to time).}
    \label{fig:cpu_time_svd}
\end{figure}

The computational cost of POD increases with the dimension $d$ of the parameter space $\caldf\subset \mathbb{R}^d$ due to the growth of the snapshot matrix $\calsf_\Phi$. 
If we collect $N^t$ time snapshots per set of parameter and require $M$ parameter instances to cover the $d$-dimensional parameter space, the total number of snapshots is $N^s = N^t\cdot M$. Since $M$ grows exponentially with $d$, the number of columns of $\calsf_\Phi$ increases accordingly. Note that the time complexity of applying the SVD to $\calsf_\Phi$ is $O(N_CN^sK)$ which leads to a substantial computational cost for high-dimensional parameter spaces. 
For example, Fig.~\ref{fig:cpu_time_svd} shows the computational time (in seconds) required to apply POD to a snapshot matrix for problem \eqref{eq:fil1}-\eqref{eq:fil6}, with $N_C = 8192$ and $N^s\in\{401, 3609, 18405, 90225\}$ 
for the dimension of the physical parameter space $d$ increasing from 0, indicating absence of physical parameters (i.e., time is the only parameter), to 3, including $\delta$, $\sigma$, and $Fr$, in addition to time. 
As expected, we observe exponential growth in the computational time as the dimension $d$ of the parameter space is increased. 

One way to reduce this growing computational cost 
is to resort to  \textit{randomized singular value decomposition} (rSVD) \cite{Martinsson2016}. Randomized POD (rPOD), i.e., the POD algorithm that replaces SVD with rSVD, is described next.

\subsection{Randomized Proper Orthogonal Decomposition} \label{sec:rpod}


As explained in the previous section, to form the reduced basis we need
only a few columns (i.e., $N_\Phi^r\ll K$) of $\caluf_\Phi$ for 
$\Phi\in\{q_l',\psi_l'\}, l=1,2$. With the rPOD, 
we wish to reconstruct the first $N_\Phi^r$ columns of the matrix $\caluf_\Phi$ efficiently and accurately.

First, we construct a Gaussian sketching matrix $\Lambda\in\mathbb{R}^{N^s\times N_\Phi^r}$, 
whose entries are drawn from a Gaussian distribution with mean $0$ and variance $1$. Matrix $\Lambda$ allows us to compress the snapshot matrix $\calsf_\Phi$, which will in turn reduce the computational cost of subsequent operations. The sketched (or sampled) matrix $\calpf_\Phi$ is computed as
\begin{equation*}
    \calpf_\Phi = \calsf_\Phi\Lambda \in \mathbb{R}^{N_C\times N_\Phi^r}.
\end{equation*}
To ensure a good reconstruction (in probability) of the basis functions, 
one typically oversamples $p$ columns, with $1\leq p\in\mathbb{N}$, and uses $\Lambda\in\mathbb{R}^{N_s\times(N_\Phi^r + p)}$ as the Gaussian sketching matrix. Typically, choosing $5\leq p\leq 20$ is sufficient to have a high probability of good reconstruction. These values work well for the best low-rank approximation of $\calsf_\Phi$ given by the deterministic truncated SVD in problems with fast decaying singular values. 

Next, we compute an orthonormal basis stored in the columns of matrix $Q_\Phi\in\mathbb{R}^{N_C\times (N_\Phi^r+p)}$ for the column space of $\calsf_\Phi$ using QR decomposition:
\begin{equation*}
    \calpf_\Phi = Q_\Phi R_\Phi, 
\end{equation*}
where $R_\Phi \in\mathbb{R}^{(N_\Phi^r + p)\times (N_\Phi^r + p)}$ is an upper triangular matrix. Then, we set
\begin{equation*}
    \caltf_\Phi = Q_\Phi^T\calsf_\Phi .
\end{equation*}
Since $\caltf_\Phi\in \mathbb{R}^{(N_{\Phi}^r + p)\times N^s}$ is a significantly smaller than $\calsf_\Phi$, SVD can be applied to $\caltf_\Phi$ more efficiently. 
The application of SVD to the matrix $\caltf_\Phi$ gives $\caltf_\Phi = \widehat{\caluf}_\Phi\Sigma_\Phi^r\calvf^{r,T}_\Phi$ where $\widehat{\caluf}_\Phi \in \mathbb{R}^{(N_\Phi^r + p)\times (N_\Phi^r + p)}$ and $\calvf_\Phi^r \in \mathbb{R}^{N^s \times (N_\Phi^r + p)}$ contain the left and right singular vectors of $\caltf_\Phi$, and $\Sigma_\Phi^r\in \mathbb{R}^{(N_\Phi^r + p)\times (N_\Phi^r + p)}$ is the diagonal matrix whose diagonal entries approximate the first $N_\Phi^r + p$ singular values of $\calsf_\Phi$.

Finally, we construct an approximation of the leading $N_\Phi^r + p$ columns of $\caluf_\Phi$ via
\begin{equation*}
    \caluf_\Phi^r = Q_\Phi \widehat{\caluf}_\Phi.
\end{equation*}
The first $N_\Phi^r$ columns of $\caluf_\Phi^r$ serves as our approximate reduced basis for variable $\Phi$.

To obtain a more accurate reconstruction of the reduced basis, we opt to use subspace iterations. This method boosts the approximation accuracy of the dominant singular vectors by amplifying the influence of the largest singular values \cite{Gu2015}. The complete algorithm of rPOD with $q$ subspace iterations is given in Alg.~\ref{alg:rpod}. 



\begin{algorithm}
    \caption{randomized POD with subspace iterations}
    \label{alg:rpod}
    \hspace*{\algorithmicindent} \textbf{Input:} snapshot matrix $\calsf_\Phi$, target rank $N_\Phi^r$, oversampling  columns $p\geq0$, and number of power iterations $q\geq1$
    \begin{algorithmic}[1]
        \State Draw a Gaussian sketching matrix $\Lambda \in \mathbb{R}^{N^s\times(N_\Phi^r + p)}$ with mean $0$ and variance $1$.
        \State Compute the sampled matrix: $\calpf_\Phi = \calsf_\Phi\Lambda$
        \State Perform $q$ subspace iterations: $\calpf_\Phi \leftarrow (S_\Phi S_\Phi^T)^q \calpf_\Phi$
        \State Apply QR decomposition: $\calpf_\Phi = Q_\Phi R_\Phi$
        \State Set $\caltf_\Phi = Q_\Phi^T S_\Phi$.
        \State Compute the SVD of $\caltf_\Phi$: $\caltf_\Phi = \widehat{\caluf}_\Phi\calsigf_\Phi^r\calvf_\Phi^{r,T}$
        \State Set $\caluf_\Phi^r = Q_\Phi\widehat{\caluf}_\Phi$.
    \end{algorithmic}
     \hspace*{\algorithmicindent} \textbf{Output:} $\caluf_\Phi^r \in \mathbb{R}^{N_h\times (N_\Phi^r + p)}, \calsigf_\Phi^r \in \mathbb{R}^{(N_\Phi^r + p)\times (N_\Phi^r + p)}, \calvf_\Phi^{r,T} \in \mathbb{R}^{(N_\Phi^r + p) \times N^s}$
\end{algorithm}

For a thorough discussion on this topic, other applications of rSVD and related methods, the readers are referred to \cite{HMT2011, Gu2015}.

\subsection{Long short-term memory network} \label{sec:lstm}

Once the reduced bases are generated, we need to find modal coefficients $\alpha_{l,i}$ and $\beta_{l,i}$ in \eqref{eq:q_approx}-\eqref{eq:psi_approx} to obtain the ROM
approximation. For this, we follow a strategy introduced in 
\cite{Besabe2024b}, which is briefly summarized below.
For more details on the use of LSTM for ROMs of 
the quasi-geostrophic systems, see also \cite{Rahman-2019, Golzar2023}.

The first step of the procedure is to approximate the snapshots using \eqref{eq:q_approx}-\eqref{eq:psi_approx}:
\begin{align}
    q_l(t_p,\xbf, \bm{\mu}_k)\approx q_{l}^r(t_p,\xbf, \bm{\mu}_k) &= \tildeq_l(\xbf, \bm{\mu}_k) + \sum_{i=1}^{N_{q_l}^r} \alpha_{l,i}(t_p, \bm{\mu}_k)\varphi_{l,i}(\xbf),\label{eq:q_approx_snap_mu} \\
    \psi_l(t_p,\xbf, \bm{\mu}_k)\approx \psi_{l}^r(t_p, \xbf, \bm{\mu}_k) &= \tildep_l(\xbf, \bm{\mu}_k) + \sum_{i=1}^{N_{\psi_l}^r} \beta_{l,i}(t_p, \bm{\mu}_k)\xi_{l,i}(\xbf), \label{eq:psi_approx_snap_mu}
\end{align}
for $l = 1,2$, $k = 1, \dots, M$, and $p = 1, \dots, N^t$, where
the coefficients $\alpha_{l,i}(t_p, \bm{\mu}_k)$ and $\beta_{l,i}(t_p, \bm{\mu}_k)$ for each time instant $t_p$ and parameter $\bm{\mu}_k$ are $i$th row and $j$th column entries of matrix $\calcf_\Phi = \caluf_{\Phi, N_{\Phi}^r}^T\calsf_\Phi\in\mathbb{R}^{N_\Phi^r\times N^s}$, where $j=(k-1)\cdot N^t + p$. 

With the modal coefficients in \eqref{eq:q_approx_snap_mu}-\eqref{eq:psi_approx_snap_mu}, we train a long short-term memory (LSTM) network, a type of recurrent neural network (RNN) introduced in \cite{LSTM1997}. LSTMs avoid
vanishing or exploding gradients that can occur in standard
RNN by employing gating functions and states to regulate the flow of information in the network. 
Like any RNN, LSTM contains recurrent layers and neurons. Within these layers, there are memory blocks, also called LSTM cells. 
Each cell contains an input gate, which adds input information to the cell, a forget gate, which discards some information from the memory to prevent overfitting, and an output gate, which passes these information to the next cell. For more information on this, we refer the reader to \cite{Besabe2024b}.

Other than the number of layers and cells per layer, one key hyperparameter of a LSTM is the lookback window $\sigma_L$, which indicates how many steps in the history of the data series the LSTM considers, both in training and prediction.
To specify the input to the network, 
let us consider the reduced space for the potential vorticity in 
layer $l$ as representative.
For simplicity of notation, we denote with $N$
the number of retained modes $N_{q_l}^r$. 
For a given time $t_p$ and parameter vector $\bm{\mu}_k$, 
we store the time information $t_p,\dots, t_{p-\sigma_L+1}$, parameter vector $\bm{\mu}_k$, and coefficients 
$\alpha_{l,i}(t_p, \bm{\mu}_k)$ in \eqref{eq:q_approx_snap_mu}
in a $\sigma_L\times (N + d + 1)$ matrix
\begin{equation}\label{eq:input}
    \left[\begin{array}{ccccc}
        \bm{\mu}_k & t_p &\alpha_{l,1}(t_p, \bm{\mu}_k) & \cdots & \alpha_{l,N}(t_p, \bm{\mu}_k) \\
        \vdots & \vdots & \vdots & \ddots & \vdots \\
        \bm{\mu}_k & t_{p-\sigma_L+1, \bm{\mu}_k} & \alpha_{l,1}(t_{p-\sigma_L + 1}, \bm{\mu}_k) & \cdots & \alpha_{l,N}(t_{p-\sigma_L + 1}, \bm{\mu}_k)
    \end{array}\right].
\end{equation}
Through the LSTM architecture, we aim to match the
training matrix \eqref{eq:input}, i.e., the
input at time $t_p$ and for parameter vector $\bm{\mu}_k$, to the corresponding output training vector $[\alpha_{l,1} (t_{p+1}, \bm{\mu}_k),\dots,$ $\alpha_{l,K}(t_{p+1}, \bm{\mu}_k)]$,
i.e., the coefficients for the reduced
order approximation at time $t_{p+1}$. 

We will denote with $\calmf_{q_l}$ the LSTM network corresponding to the potential vorticity $q_l$, while the LSTM network associated with the stream function $\psi_l$ is denoted by $\calmf_{\psi_l}$.

\section{Numerical Results} \label{sec:num_res}

For our numerical tests, we consider an extension \cite{Besabe2024, Besabe2024b} of the so-called double-wind gyre forcing experiment, which is a classical benchmark for new numerical models of geophysical flows \cite{Nadiga2001, Holm2003, Greathbatch2000, Monteiro2014, Monteiro2015, Girfoglio_JCAM2023, QGE-review, San2015, Girfoglio2023}. The computational domain is $\Omega = [0,1]\times[-1,1]$ and the wind forcing is given by $F=\sin(\pi y)$. We set $Re = 450$ and $Ro = 0.001$.
For the full order simulations, we use a mesh with size $h=1/64$, which is slightly larger than the Munk scale,  and set the filtering radius to $\alpha_1 = \alpha_2 = h$.
The time step is set to $\Delta t = 2.5E-05$. From the full order simulations, we collect 401 snapshots over the training time interval $[10,50]$, i.e., 
we sample each variable every $0.1$ time unit.

In addition to time,  
we have three physical variable parameters: aspect ratio $\delta$, 
friction coefficient on the ocean floor $\sigma$, 
and Froude number $Fr$. 
For the physical parameters, we will consider the following ranges:
\begin{align*}
    \delta \in [0.2, 0.6], \quad \sigma\in [0.006, 0.01], \quad Fr\in [0.07, 0.11],
\end{align*}
which will be motivated in Sec.~\ref{sec:time+delta}-\ref{sec:time+delta+sigma+Fr}. We note that these physical parameters
are different in nature: one characterizes the geometry ($\delta$), another the physics of the flow ($Fr$), and
the last is a non-dimensionalized physical parameter ($\sigma$). In addition, they have different magnitudes. 
So, we will gradually increase the dimension of the parameter space from 2 
(time and $\delta$) in Sec.~\ref{sec:time+delta}  to 3 (time, $\delta$, and $\sigma$) in Sec.~\ref{sec:time+delta+sigma} and finally
4 (time and all three physical parameters) in Sec.~\ref{sec:time+delta+sigma+Fr}. However, before starting the parametric study, 
we need to understand the proper setting of oversampling 
$p$ in the rPOD algorithm.

\begin{figure}[htb!]
\centering
 \begin{subfigure}{0.49\textwidth}
     \includegraphics[width=0.9\textwidth]{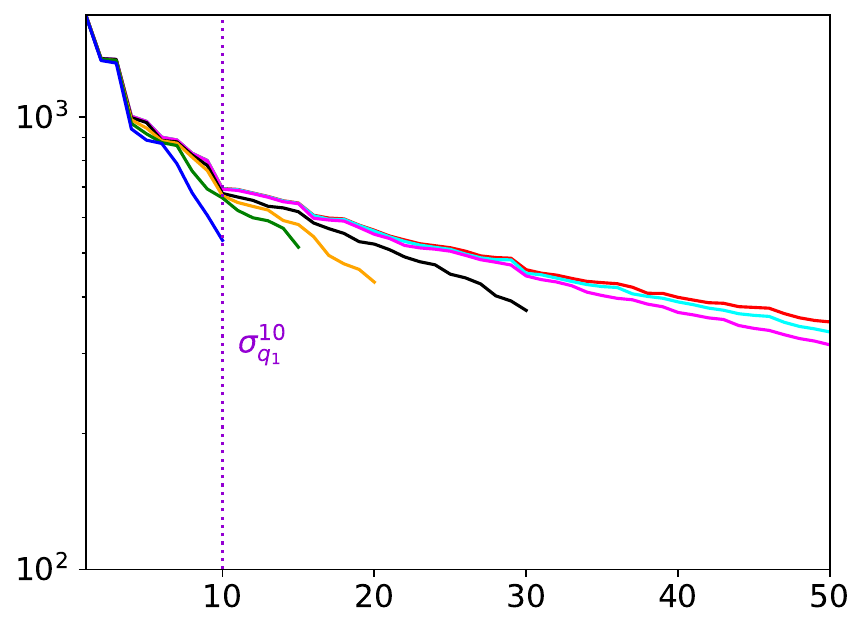}
     \label{fig:q1_sv}
     \caption{singular values for ${q_1}$}
 \end{subfigure}
 \hfill
 \begin{subfigure}{0.49\textwidth}
     \includegraphics[width=0.9\textwidth]{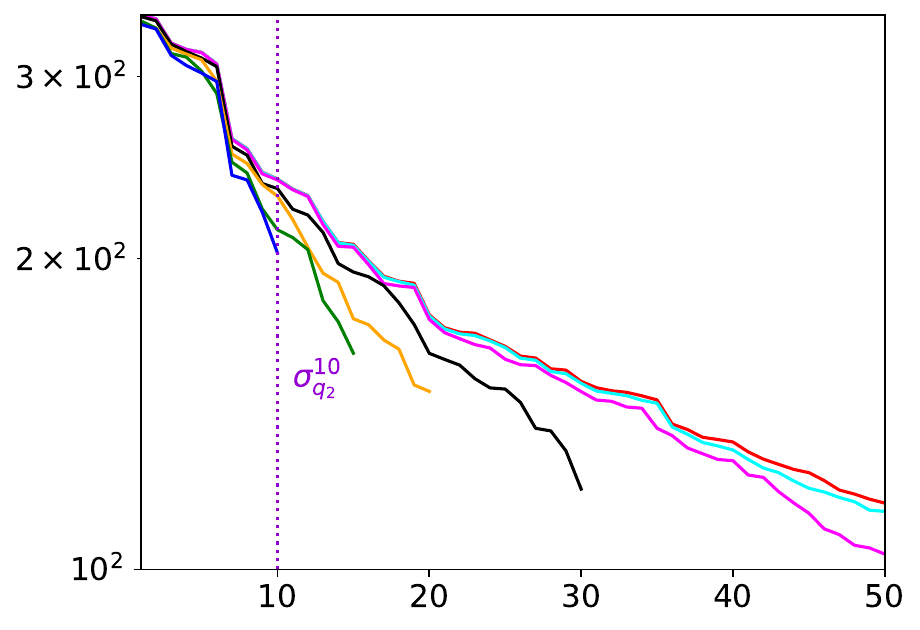}
     \label{fig:q2_sv}
     \caption{singular values for ${q_2}$}
 \end{subfigure}
\hfill
 \begin{subfigure}{0.49\textwidth}
     \includegraphics[width=0.9\textwidth]{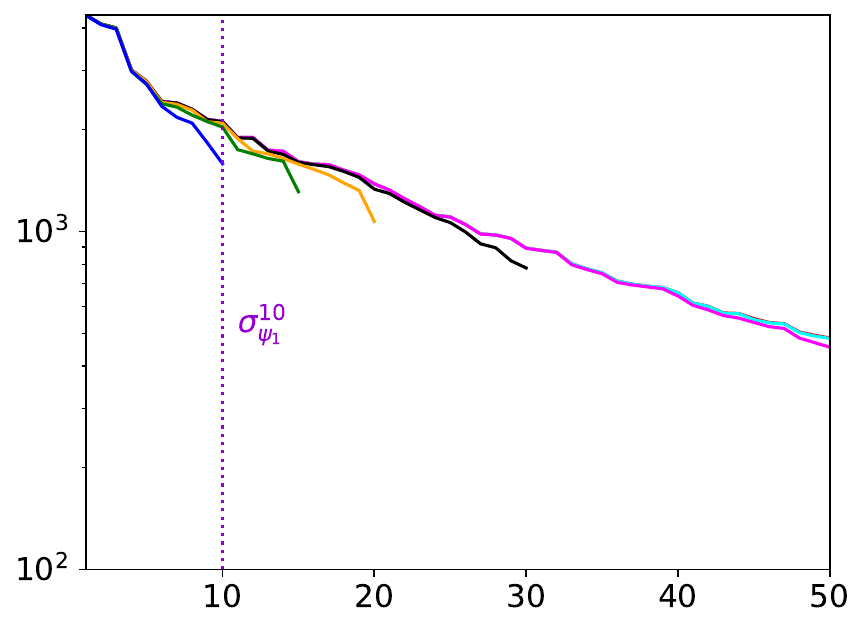}
     \label{fig:psi1_sv}
     \caption{singular values for ${\psi_1}$}
 \end{subfigure}
 \hfill
 \begin{subfigure}{0.49\textwidth}
     \includegraphics[width=0.9\textwidth]{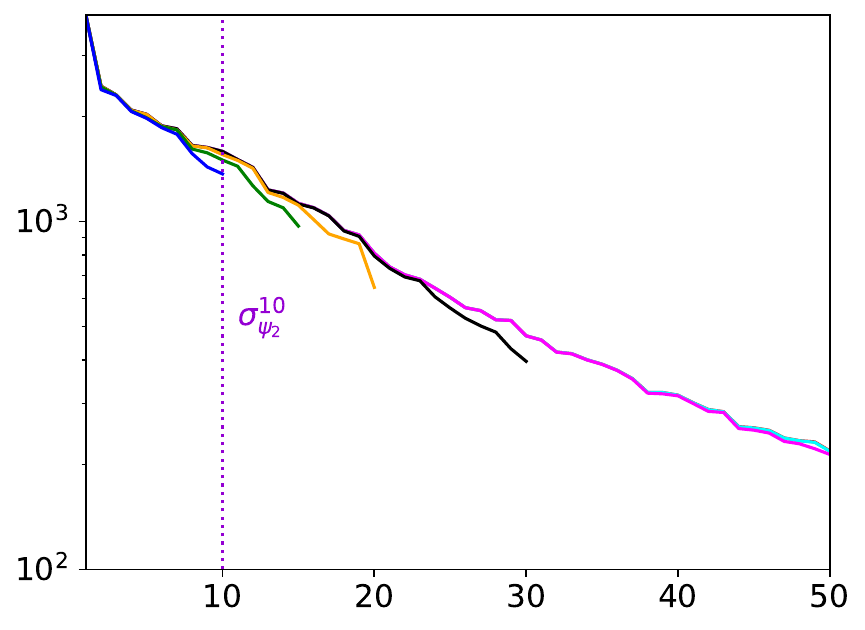}
     \label{fig:psi2_sv}
     \caption{singular values for ${\psi_2}$}
 \end{subfigure}
\vskip .2cm
\begin{subfigure}{0.95\textwidth}
    \centering
    \includegraphics[width=\textwidth]{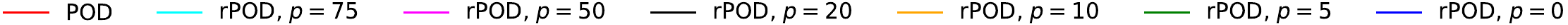}
\end{subfigure}
 \caption{
 Decay of singular values from rPOD with target rank $N_\Phi^r = 10$ and 
 $p=0,5,10,20,50,75$ for variables $q_1$ (top left), $q_2$ (top right),$\psi_1$ (bottom left) and $\psi_2$ (bottom right). 
 }
 \label{fig:singular_values}
\end{figure}

\begin{figure}
    \centering
        \begin{tabular}{ccccccccc}
             & $\varphi_{l,1}^\text{POD}$ & \hspace{-0.5cm} $\varphi_{l,1}^0$ & \hspace{-0.5cm} $\varphi_{l,1}^{75}$ & \hspace{-1.42cm} \footnotesize{$\abs{\varphi_{l,1}^{\text{POD}} - \varphi_{l,1}^{75}}$} & $\varphi_{l,10}^\text{POD}$ & \hspace{-0.5cm} $\varphi_{l,10}^0$ & \hspace{-0.5cm} $\varphi_{l,10}^{75}$ & \hspace{-1.42cm} \footnotesize{$\abs{\varphi_{l,10}^{\text{POD}} - \varphi_{l,10}^{75}}$} \\
             
            \hspace{-0.2cm}$q_1$ &\hspace{-0.4cm} \includegraphics[align=c,scale = 0.18]{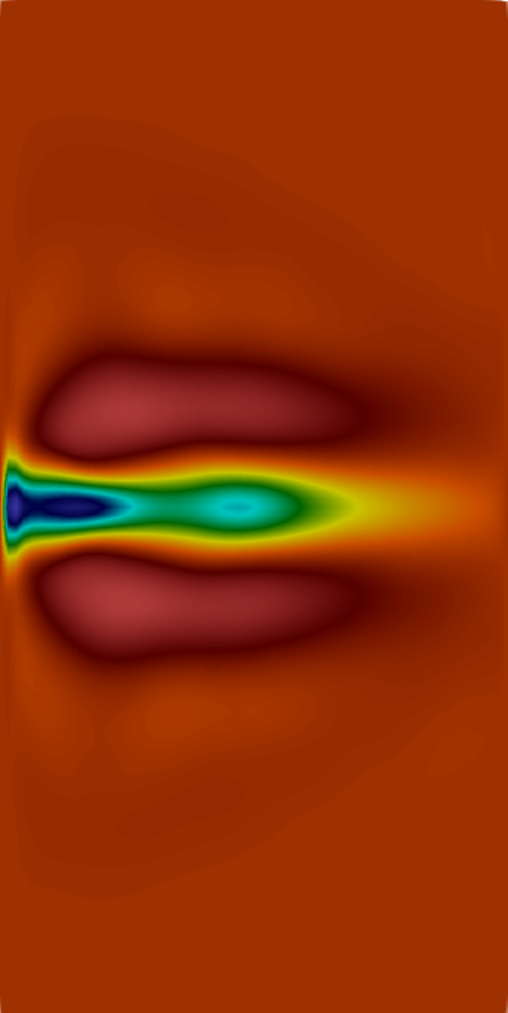} & \hspace{-0.5cm} \includegraphics[align=c,scale = 0.18]{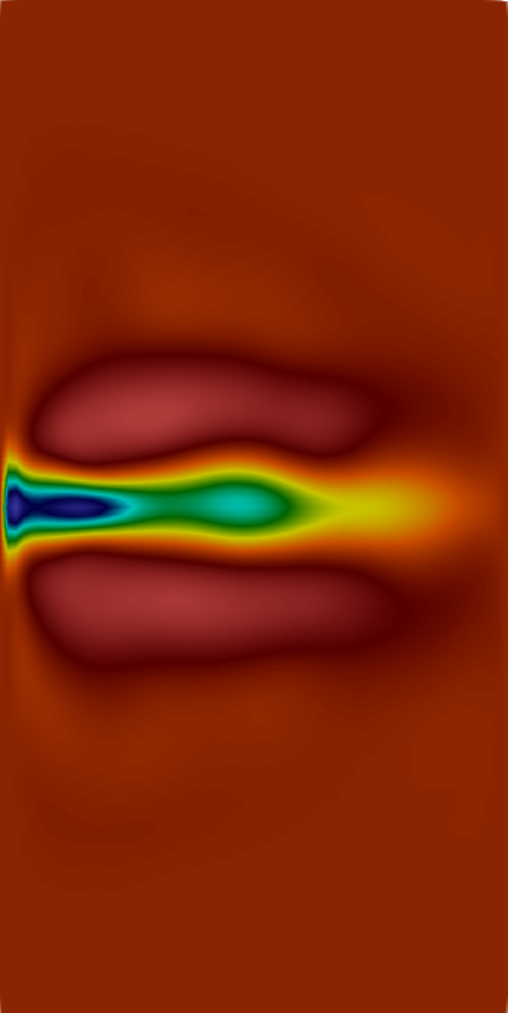} & \hspace{-0.5cm} \includegraphics[align=c,scale = 0.18]{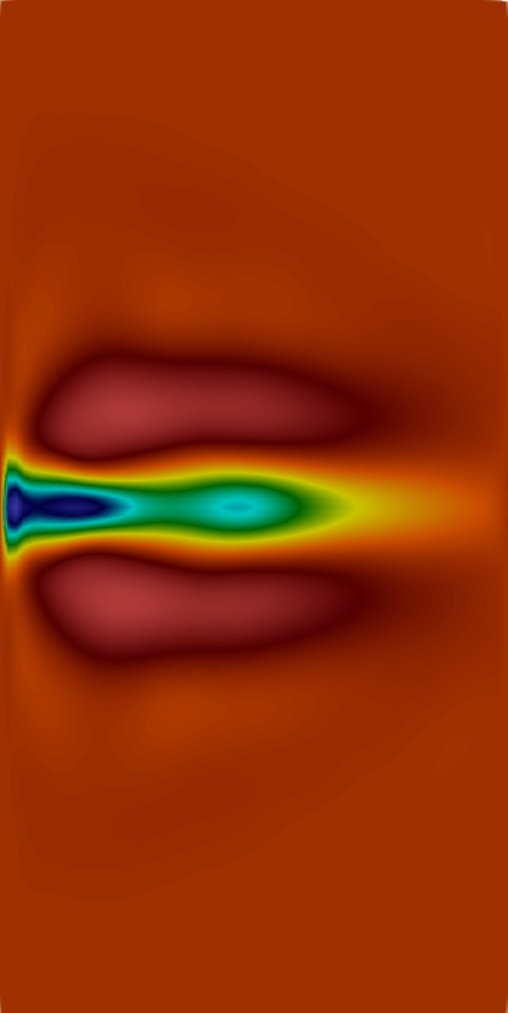} & \hspace{-0.5cm} \includegraphics[align=c,scale = 0.18]{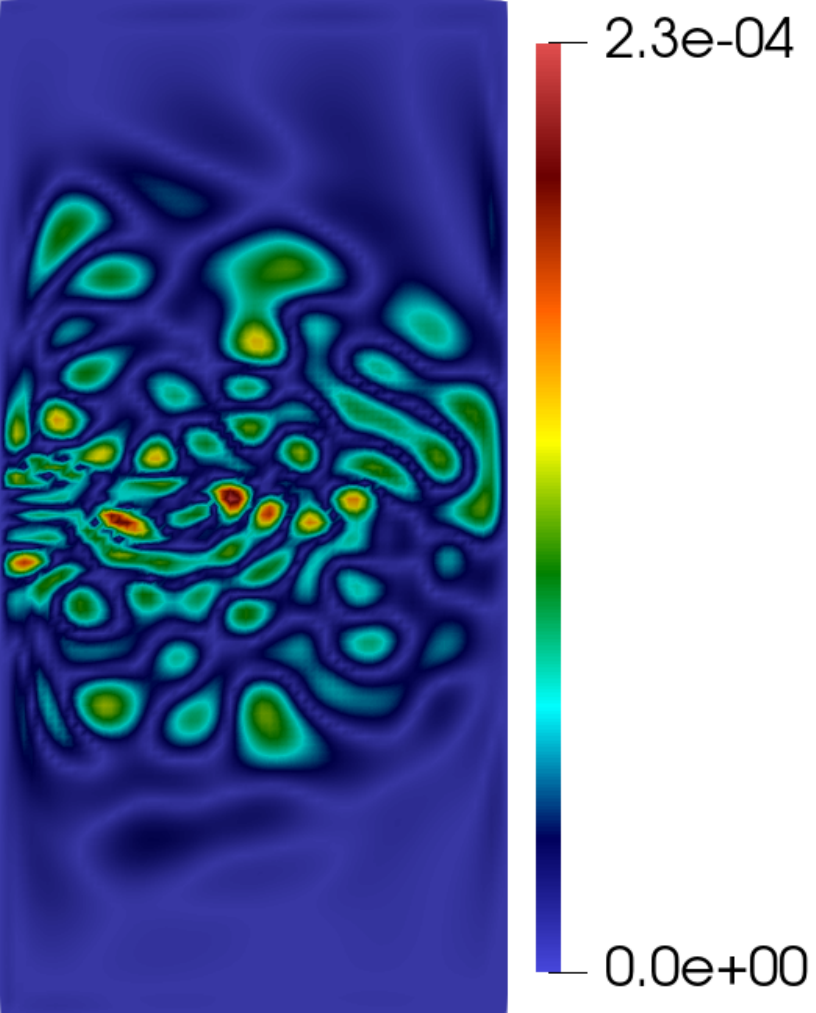} & \hspace{-0.3cm}\includegraphics[align=c,scale = 0.18]{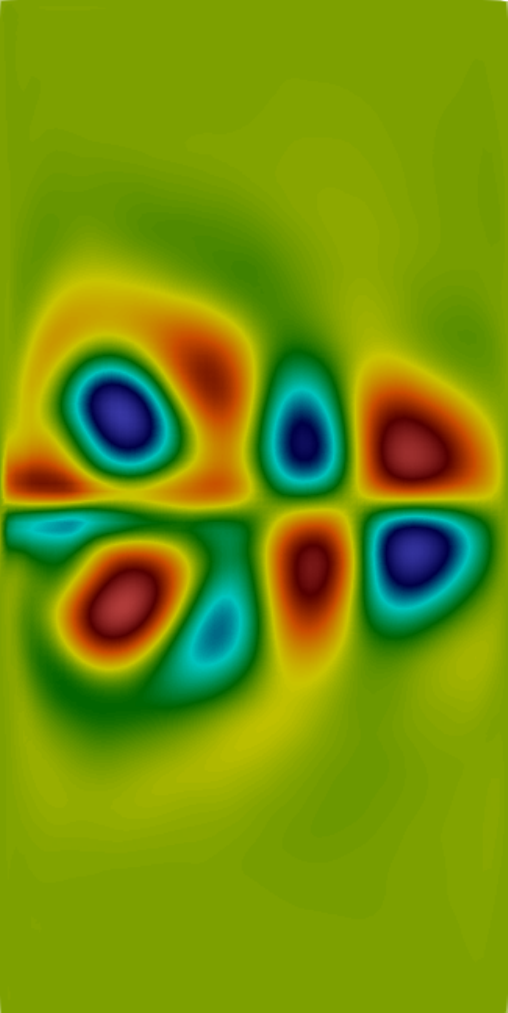} & \hspace{-0.5cm} \includegraphics[align=c,scale = 0.18]{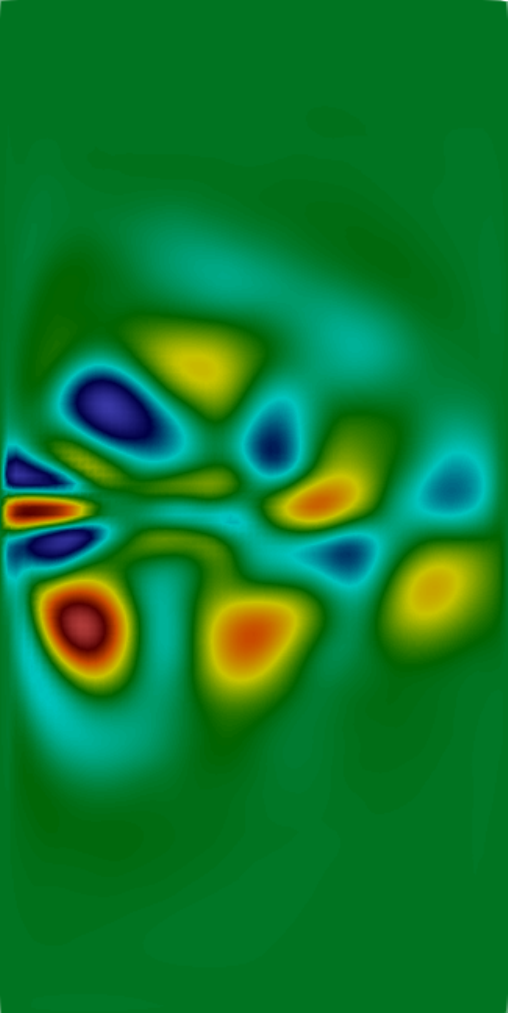} & \hspace{-0.5cm} \includegraphics[align=c,scale = 0.18]{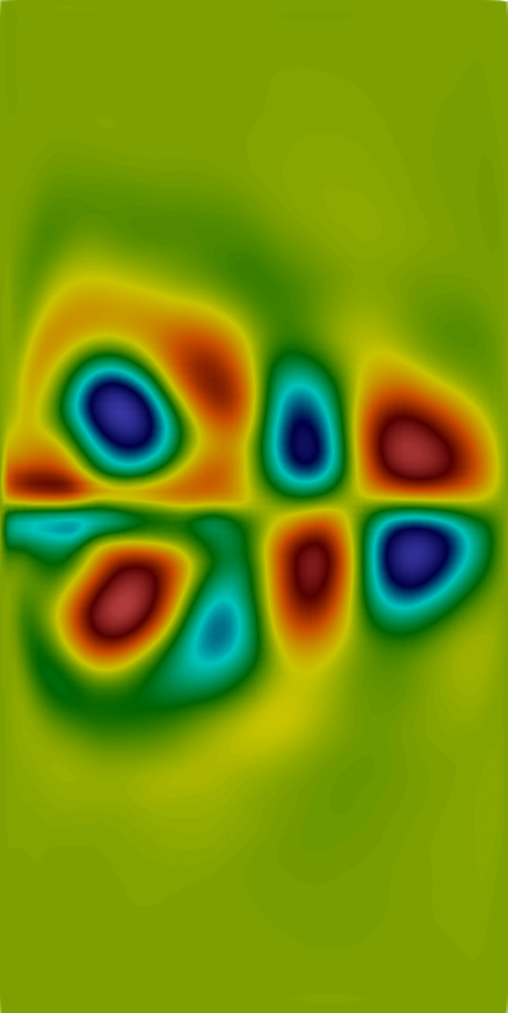} & \hspace{-0.5cm} \includegraphics[align=c,scale = 0.18]{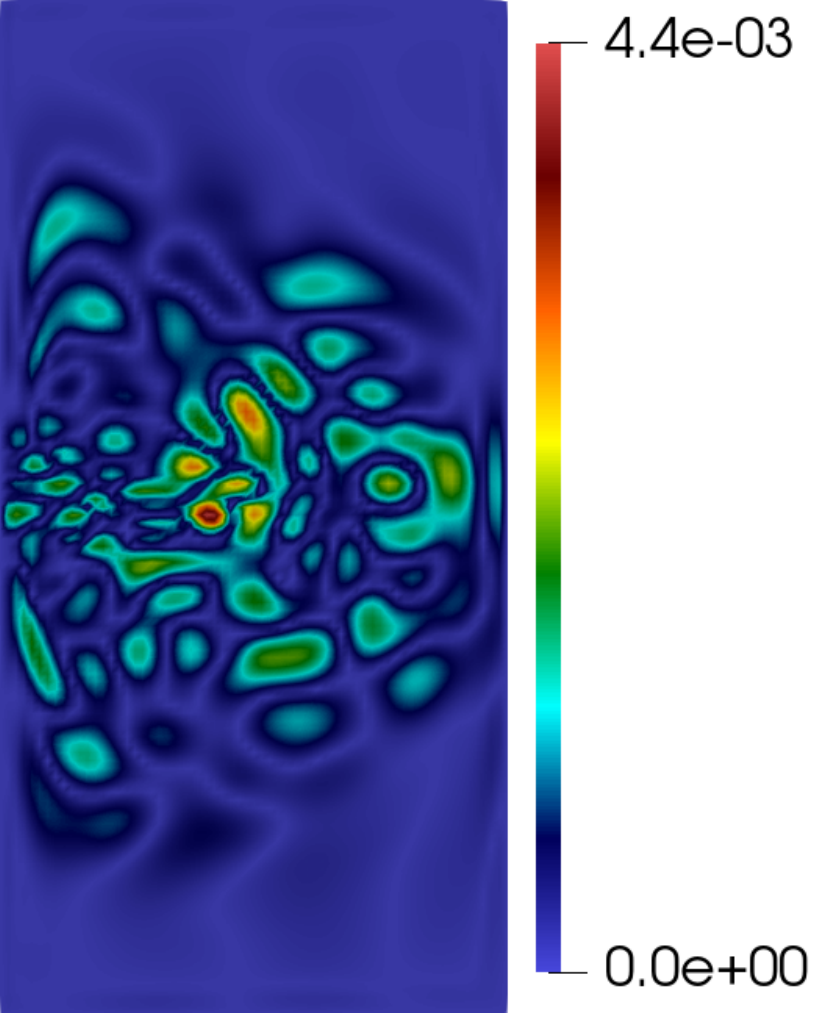} \\
            \vspace{-0.3cm}&&&&\\
            \hspace{-0.2cm}$q_2$ & \hspace{-0.3cm}\includegraphics[align=c,scale = 0.18]{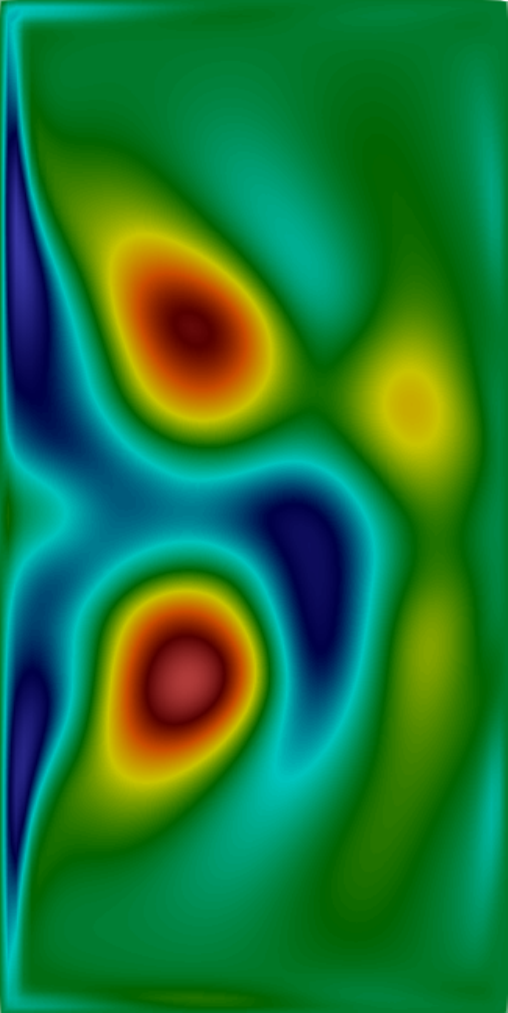} & \hspace{-0.5cm} \includegraphics[align=c,scale = 0.18]{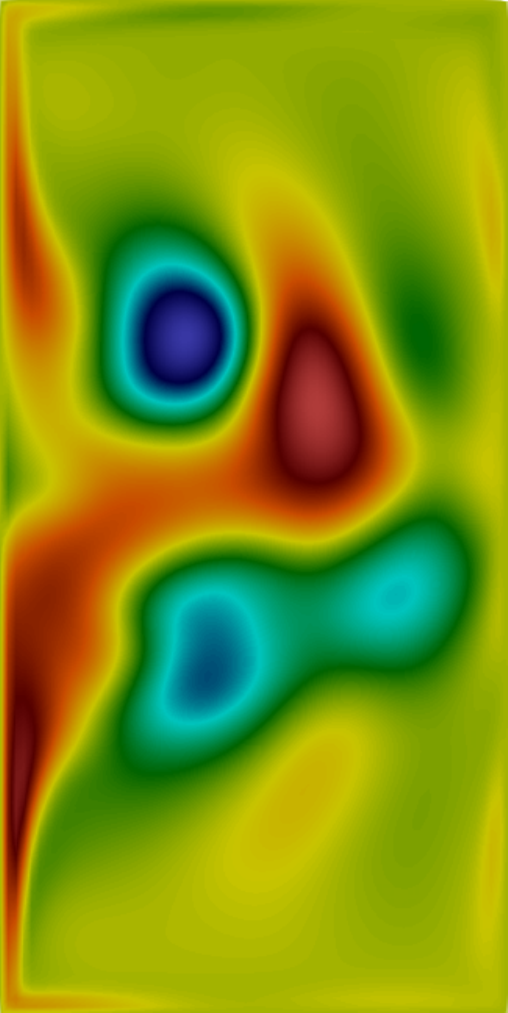} & \hspace{-0.5cm} \includegraphics[align=c,scale = 0.18]{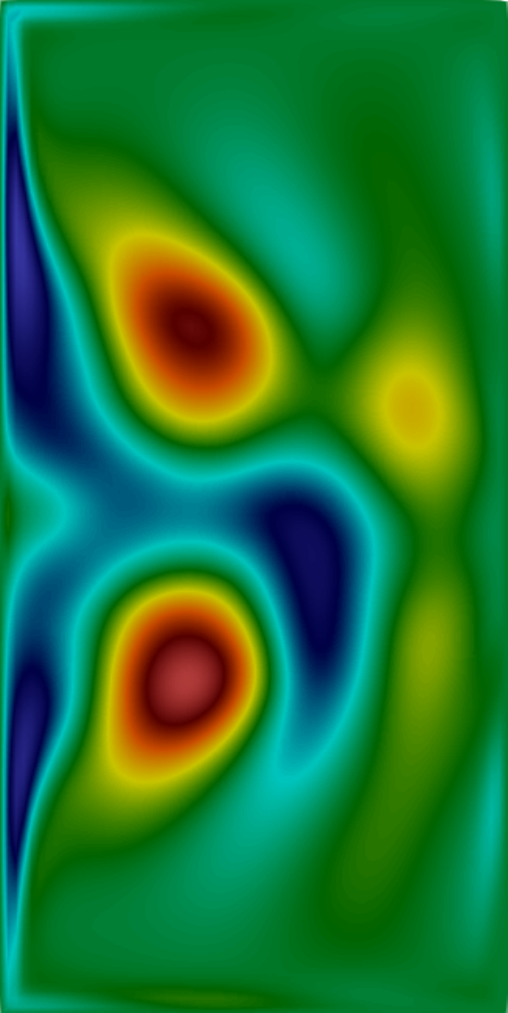} & \hspace{-0.5cm} \includegraphics[align=c,scale = 0.18]{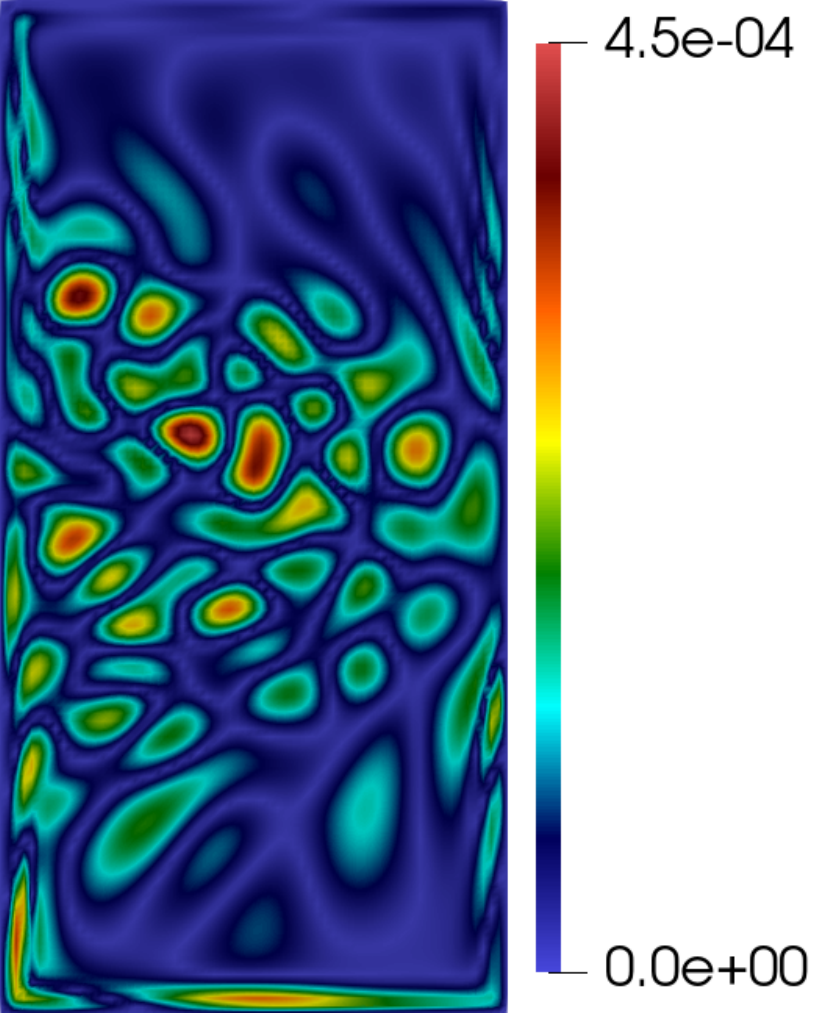} & \hspace{-0.3cm}\includegraphics[align=c,scale = 0.18]{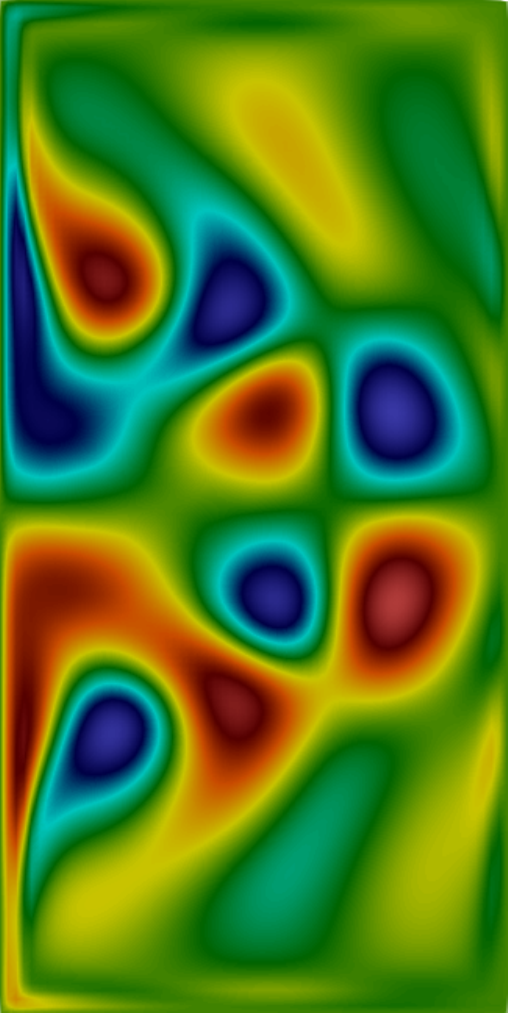} & \hspace{-0.5cm} \includegraphics[align=c,scale = 0.18]{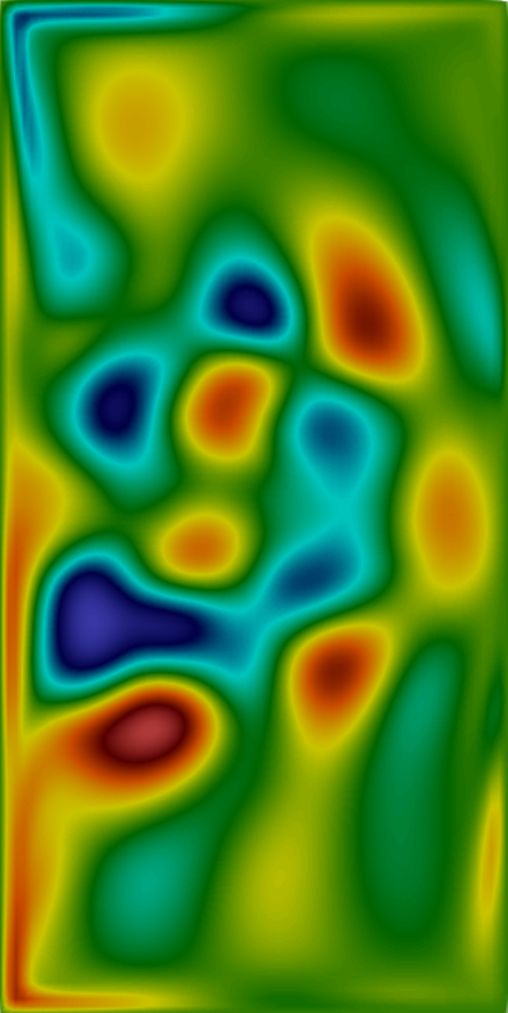} & \hspace{-0.5cm} \includegraphics[align=c,scale = 0.18]{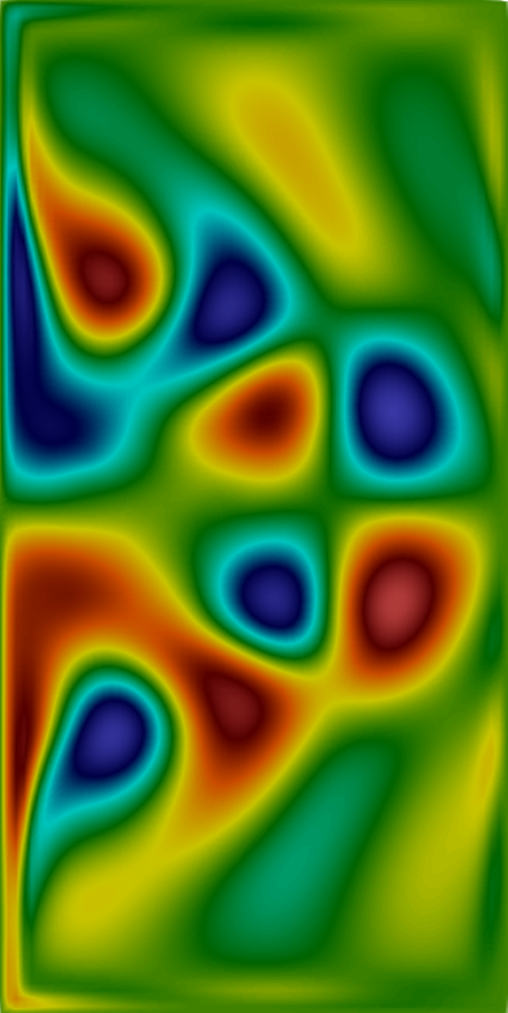} & \hspace{-0.5cm} \includegraphics[align=c,scale = 0.18]{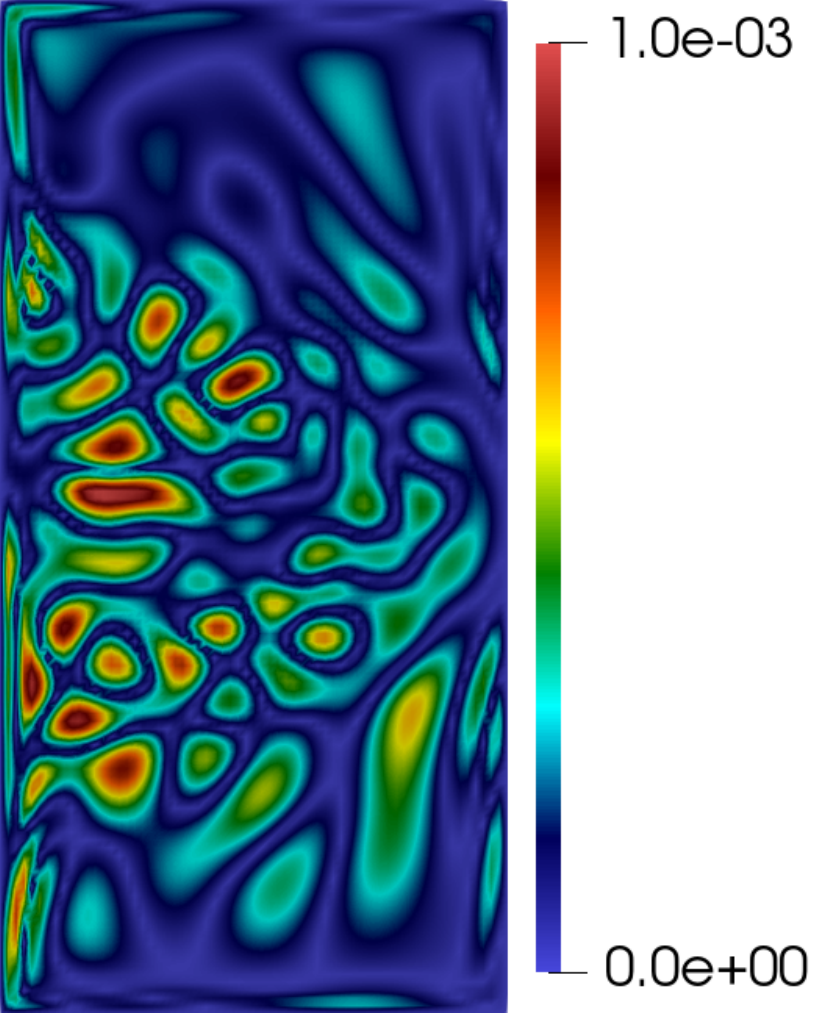} 
        \end{tabular}
    \caption{
    First (left panel) and tenth (right panel) modes of $q_1$ (first row) and $q_2$ (second row) computed using POD (first column in each panel) and rPOD with $p=0$ (second column in each panel) and oversampling of $p=75$ (third column in each panel).
    The fourth column in each panel shows the absolute difference between
    modes from POD and rPOD with $p=75$. 
    }
    \label{fig:modes_q}
\end{figure}

\begin{figure}
    \centering
        \begin{tabular}{ccccccccc}
             & $\xi_{l,1}^\text{POD}$ & \hspace{-0.5cm} $\xi_{l,1}^0$ & \hspace{-0.5cm} $\xi_{l,1}^{75}$ & \hspace{-1.42cm} \footnotesize{$\abs{\xi_{l,1}^{\text{POD}} - \xi_{l,1}^{75}}$} & $\xi_{l,10}^\text{POD}$ & \hspace{-0.5cm} $\xi_{l,10}^0$ & \hspace{-0.5cm} $\xi_{l,10}^{75}$ & \hspace{-1.42cm} \footnotesize{$\abs{\xi_{l,10}^{\text{POD}} - \xi_{l,10}^{75}}$} \\
             
            \hspace{-0.2cm}$\psi_1$ & \hspace{-0.3cm}\includegraphics[align=c,scale = 0.18]{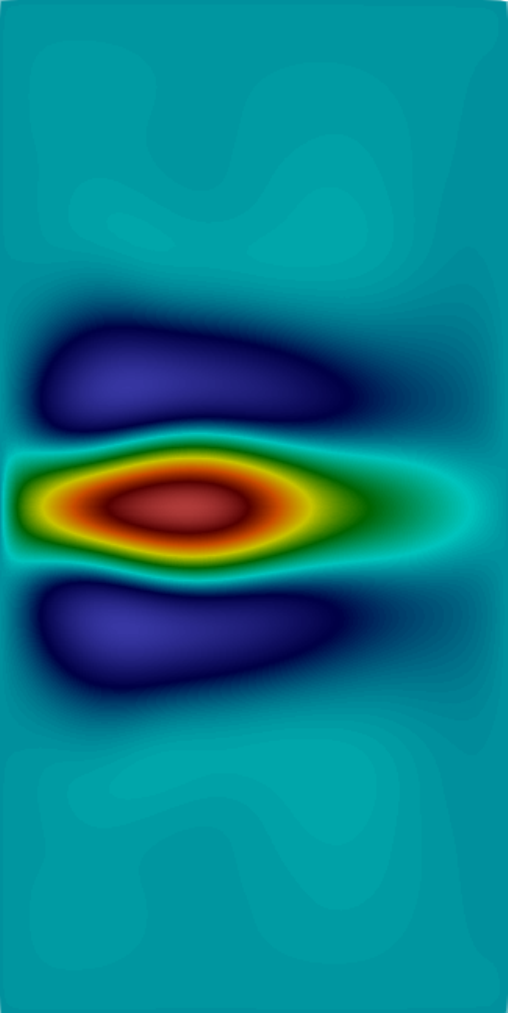} & \hspace{-0.5cm} \includegraphics[align=c,scale = 0.18]{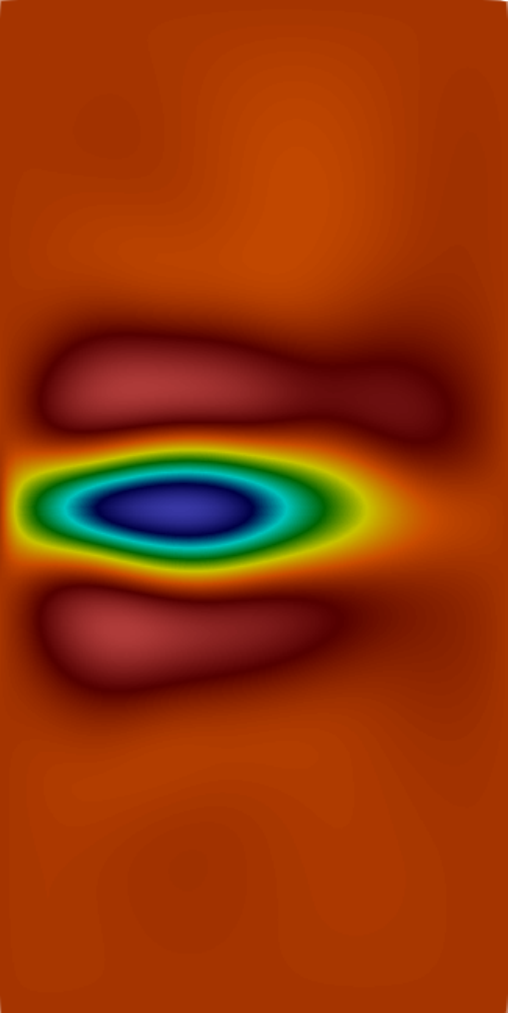} & \hspace{-0.5cm} \includegraphics[align=c,scale = 0.18]{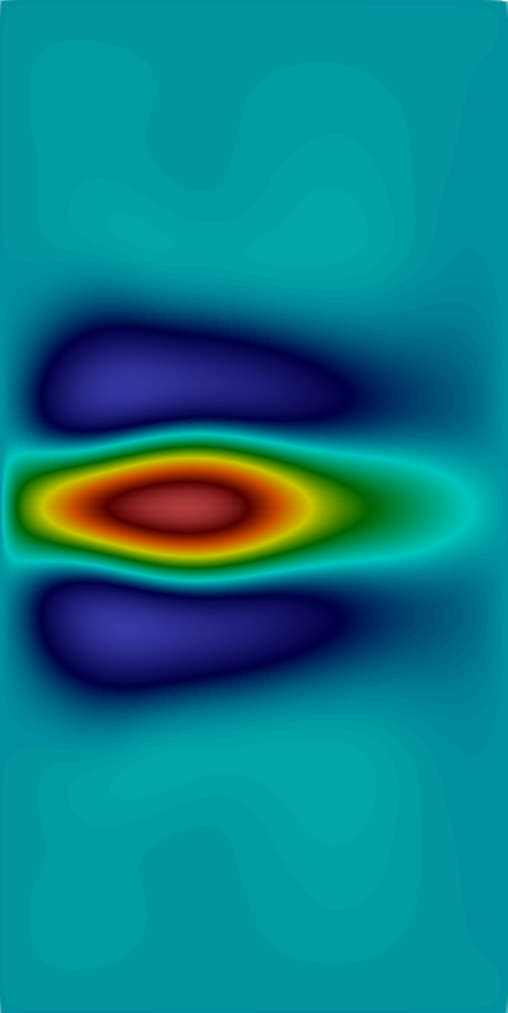} & \hspace{-0.5cm} \includegraphics[align=c,scale = 0.18]{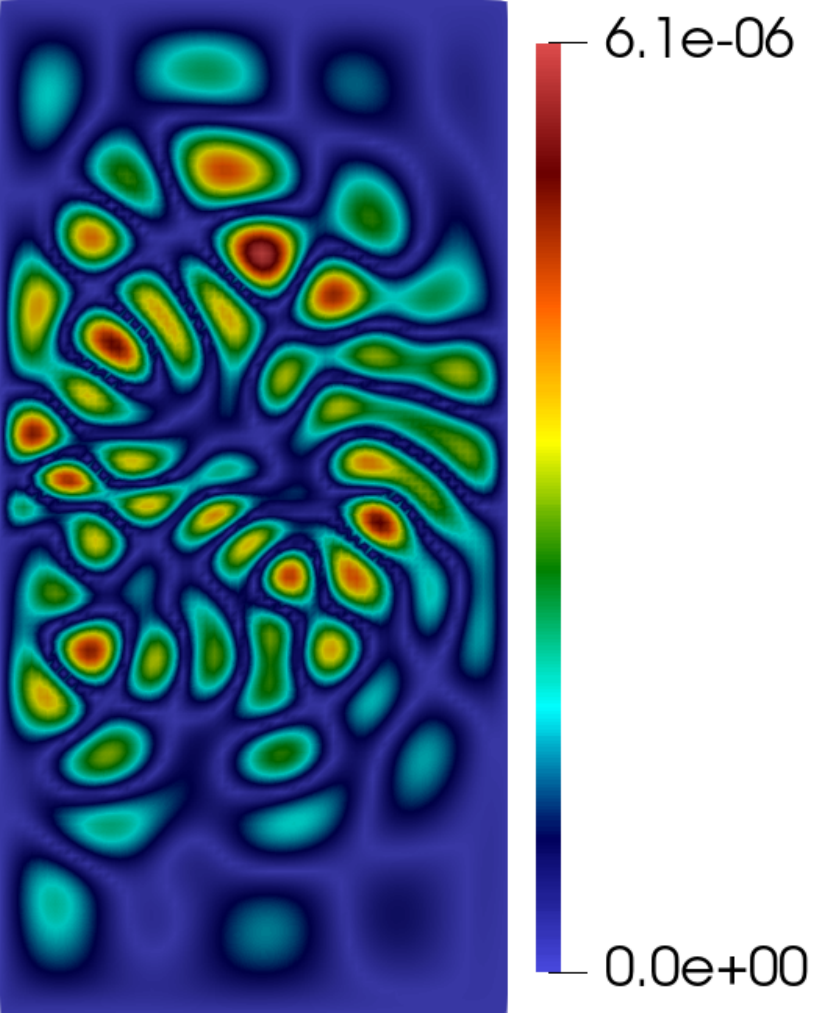} & \hspace{-0.3cm}\includegraphics[align=c,scale = 0.18]{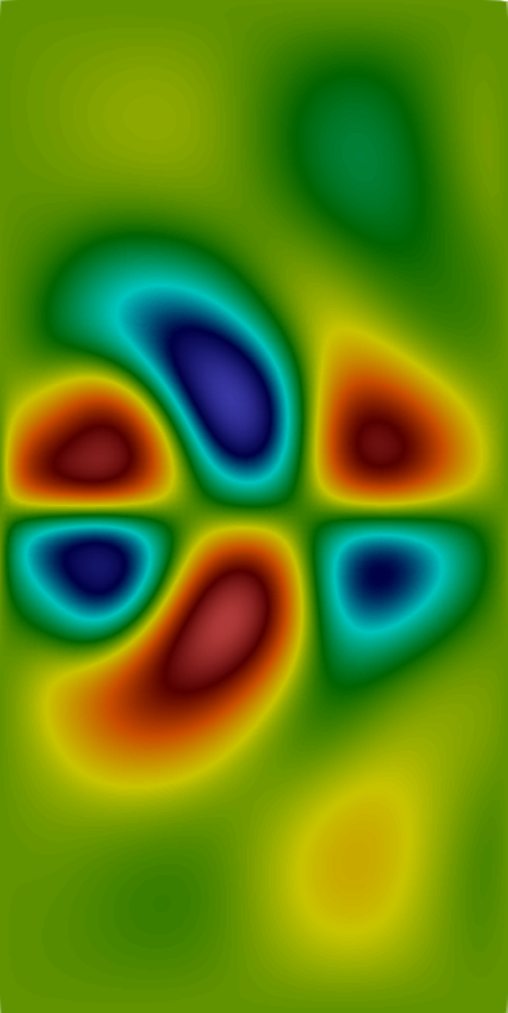} & \hspace{-0.5cm} \includegraphics[align=c,scale = 0.18]{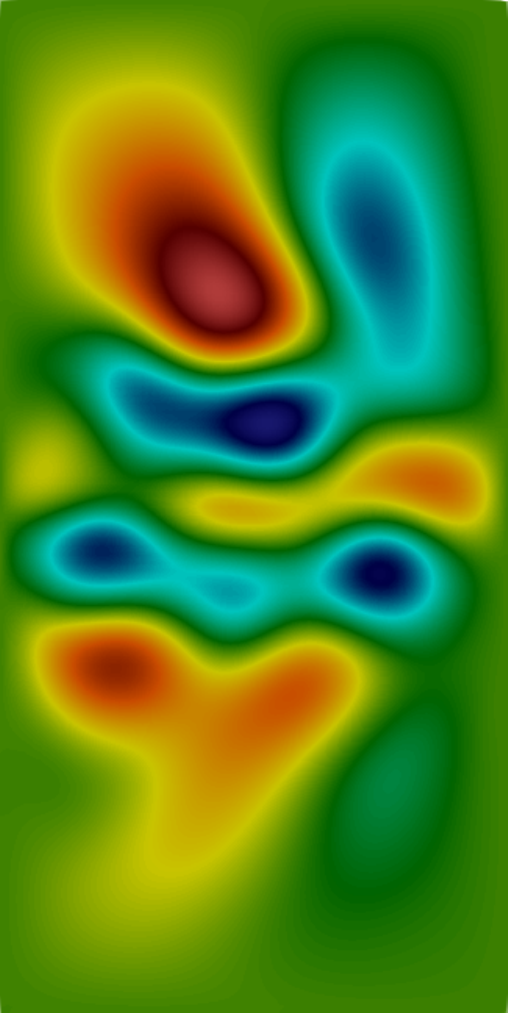} & \hspace{-0.5cm} \includegraphics[align=c,scale = 0.18]{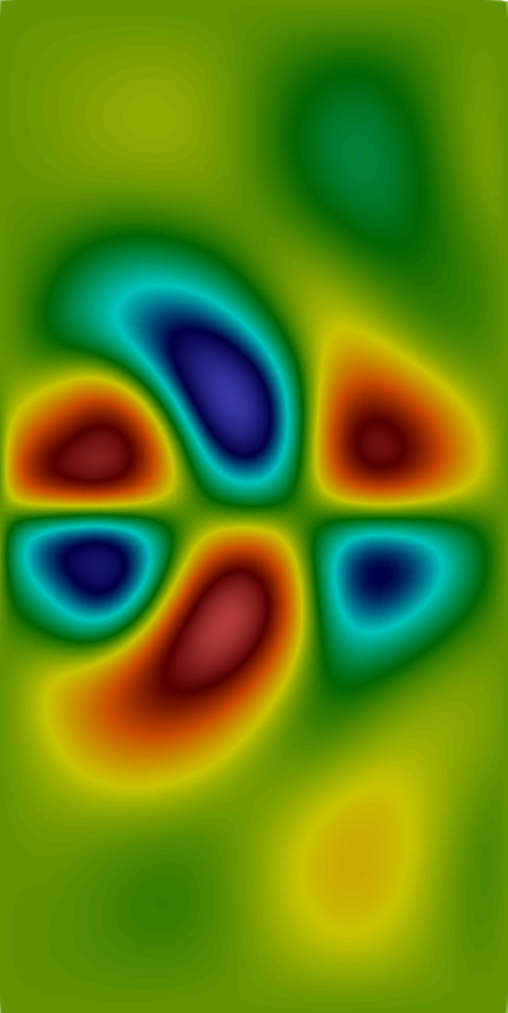} & \hspace{-0.5cm} \includegraphics[align=c,scale = 0.18]{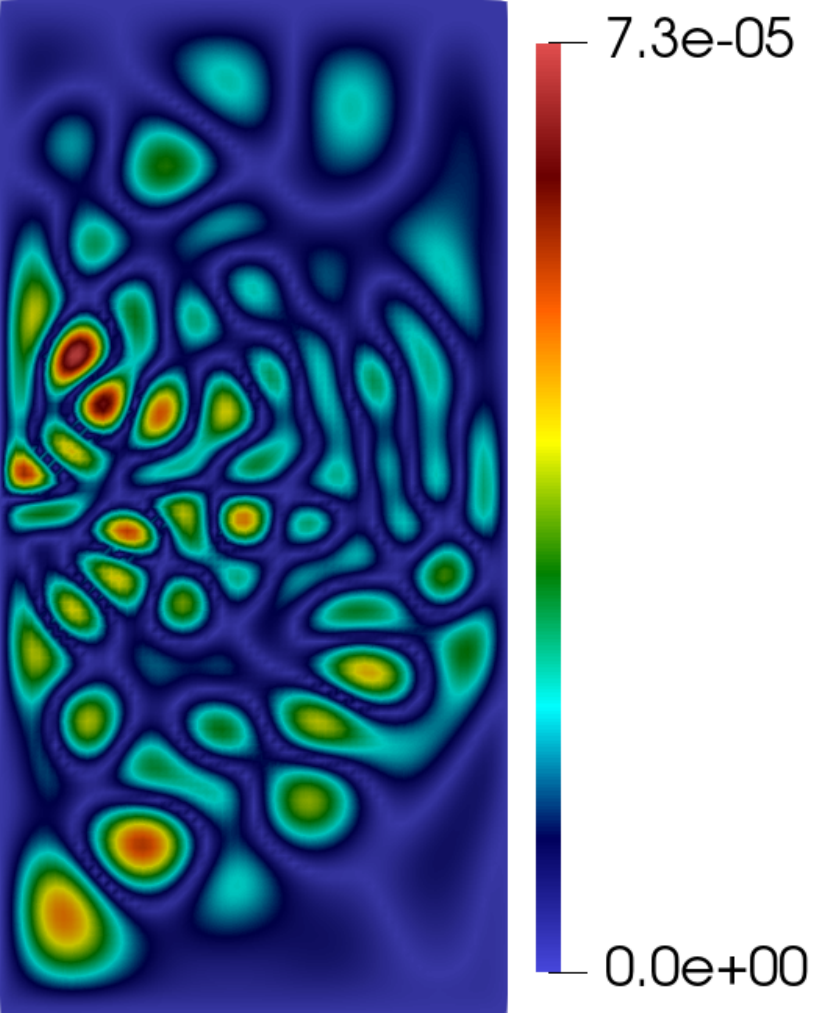} \\
            \vspace{-0.3cm}&&&&\\
            \hspace{-0.cm}$\psi_2$ & \hspace{-0.3cm}\includegraphics[align=c,scale = 0.18]{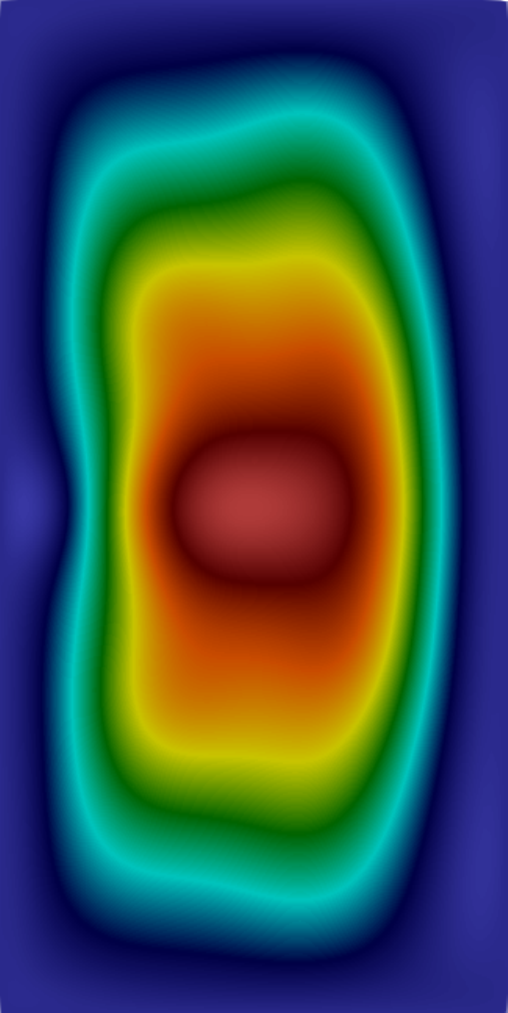} & \hspace{-0.5cm} \includegraphics[align=c,scale = 0.18]{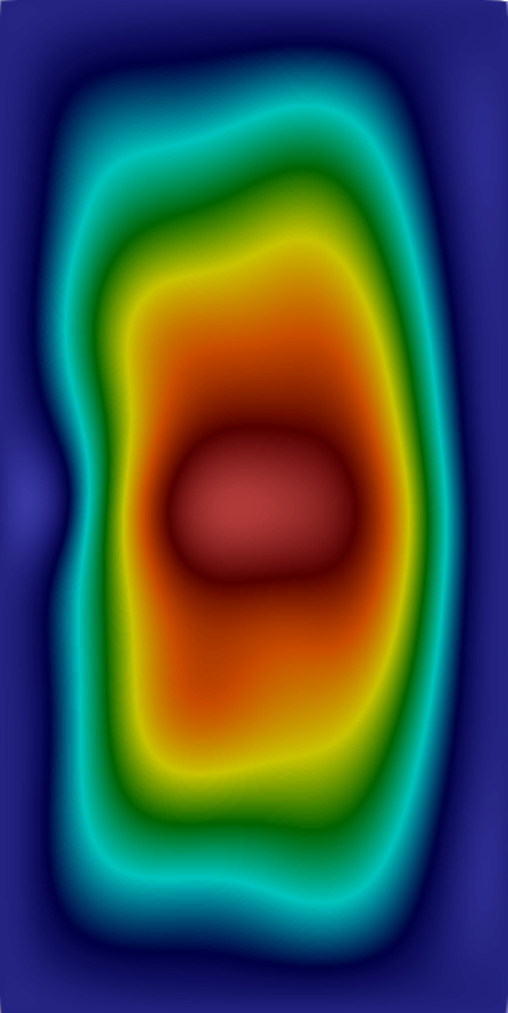} & \hspace{-0.5cm} \includegraphics[align=c,scale = 0.18]{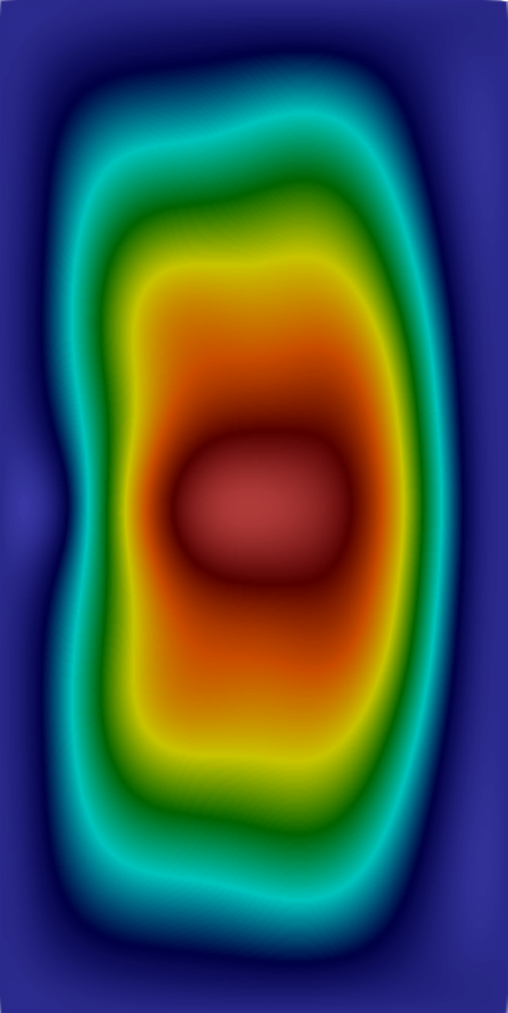} & \hspace{-0.5cm} \includegraphics[align=c,scale = 0.18]{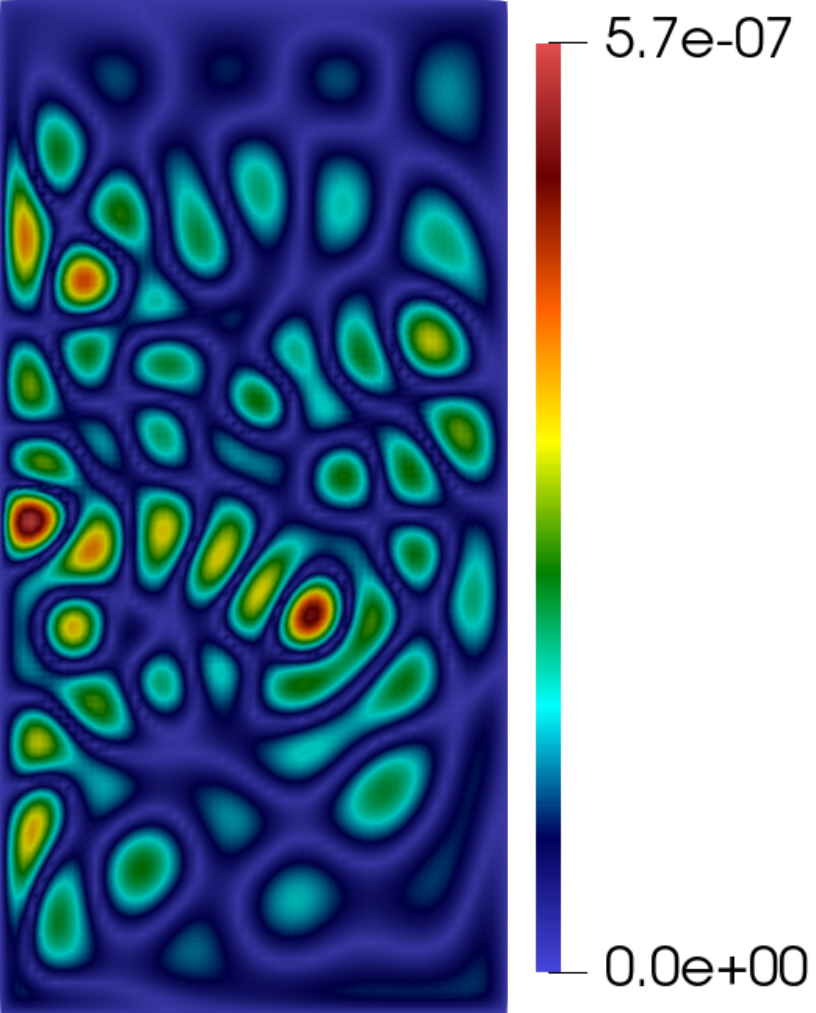} & \hspace{-0.3cm}\includegraphics[align=c,scale = 0.18]{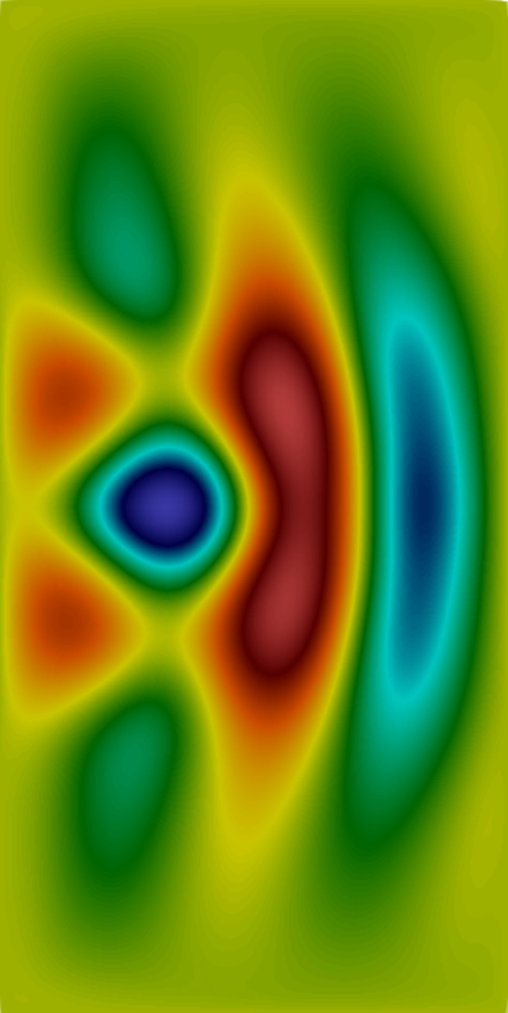} & \hspace{-0.5cm} \includegraphics[align=c,scale = 0.18]{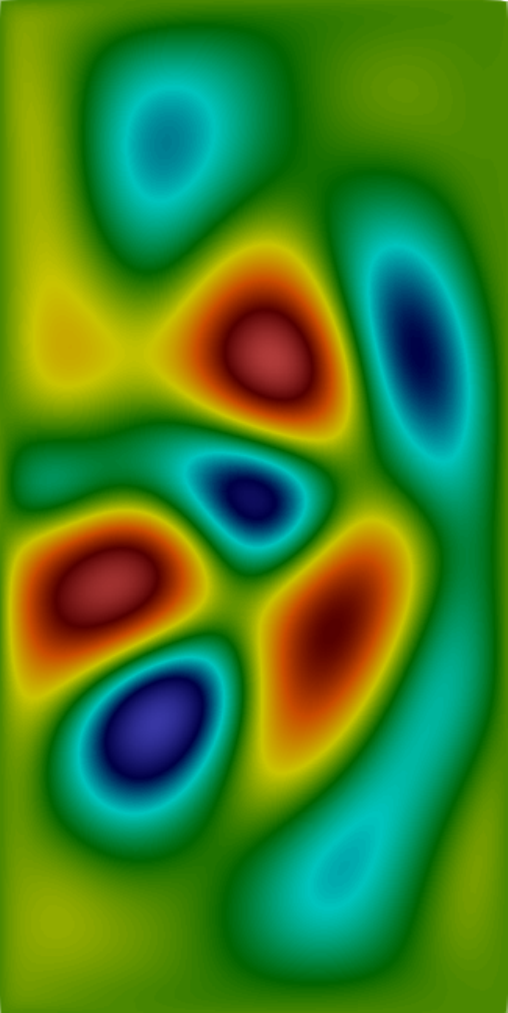} & \hspace{-0.5cm} \includegraphics[align=c,scale = 0.18]{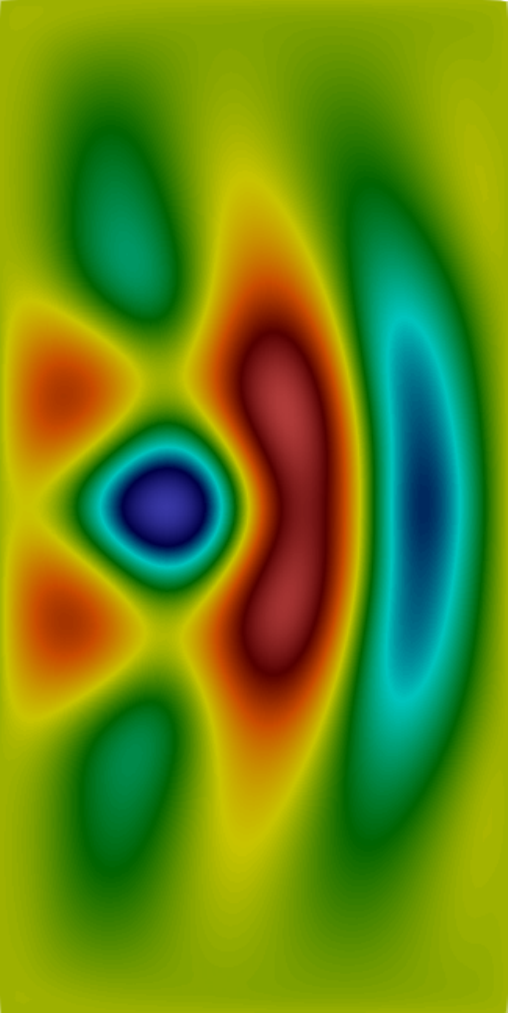} & \hspace{-0.5cm} \includegraphics[align=c,scale = 0.18]{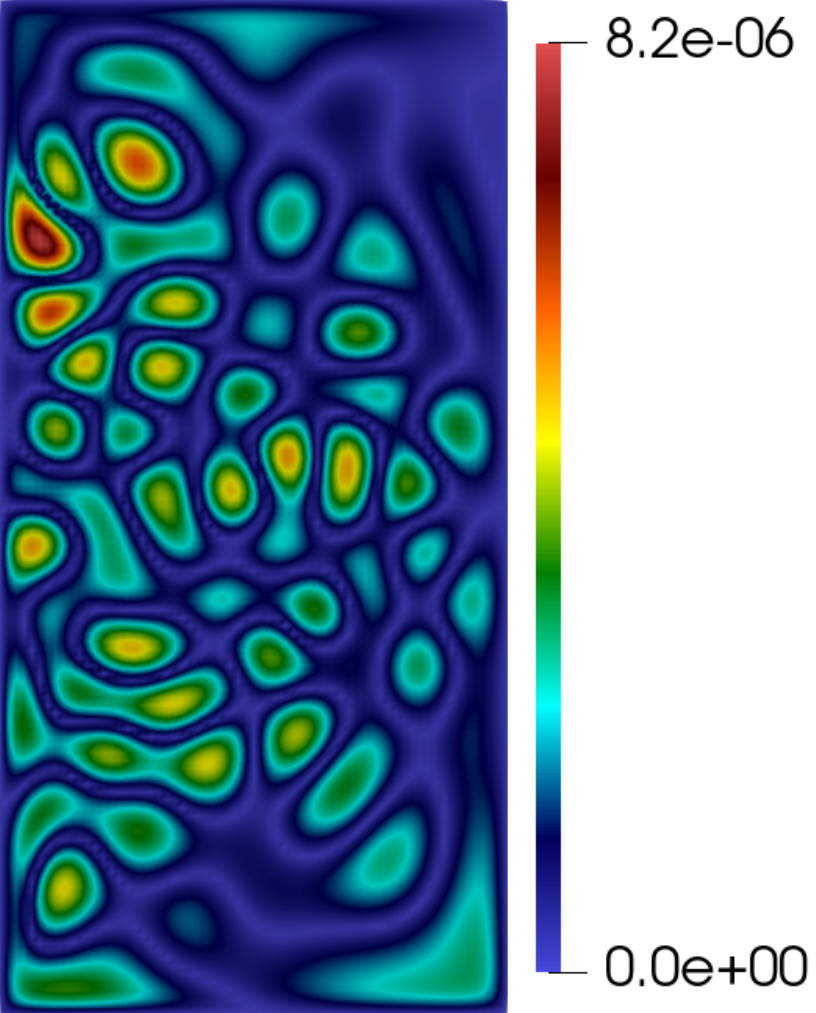} 
        \end{tabular}
    \caption{
    First (left panel) and tenth (right panel) modes of $\psi_1$ (first row) and $\psi_2$ (second row) computed using POD (first column in each panel) and rPOD with $p=0$ (second column in each panel) and oversampling of $p=75$ (third column in each panel). The fourth column in each panel shows the absolute difference between
    modes from POD and rPOD with $p=75$. 
    }
    \label{fig:modes_psi}
\end{figure}

To determine a suitable value of $p$, we consider the following sampling for the varying physical parameters:
\begin{align} \label{eq:test_delta}
    & \delta = \{0.2, 0.25, 0.3, 0.35, 0.4, 0.45, 0.5, 0.55, 0.6\}\\
    & \sigma = \{0.006, 0.007, 0.008, 0.009, 0.010\}\\
    & Fr = \{0.07, 0.08, 0.09, 0.10, 0.11\} \label{eq:test_fr}
\end{align}
Fig.~\ref{fig:singular_values} shows the decay of singular values 
computed by POD and rPOD with $q=1$ (see Alg.~\ref{alg:rpod}, l.~3) and $p=0, 5, 10, 20, 50, 75$, where $p=0$ means no oversampling. 
We see that rPOD provides a better approximation of the first 
few POD singular values as $p$ is increased. 
This is confirmed by Fig.~\ref{fig:modes_q}-\ref{fig:modes_psi}, which show a comparison of the POD and rPOD ($p = 0$ and $p = 75$) first and tenth basis functions. Note that the reconstruction of the modes is more accurate for $\psi_l$ than of $q_l$. We suspect that this is due to the faster decay of singular values associated to the stream functions $\psi_l$ (see Fig.~\ref{fig:singular_values}). Note that the vertical axis in the different
panels of Fig.~\ref{fig:singular_values} are different. 


\begin{figure}
 \begin{subfigure}{0.49\textwidth}
     \includegraphics[width=\textwidth]{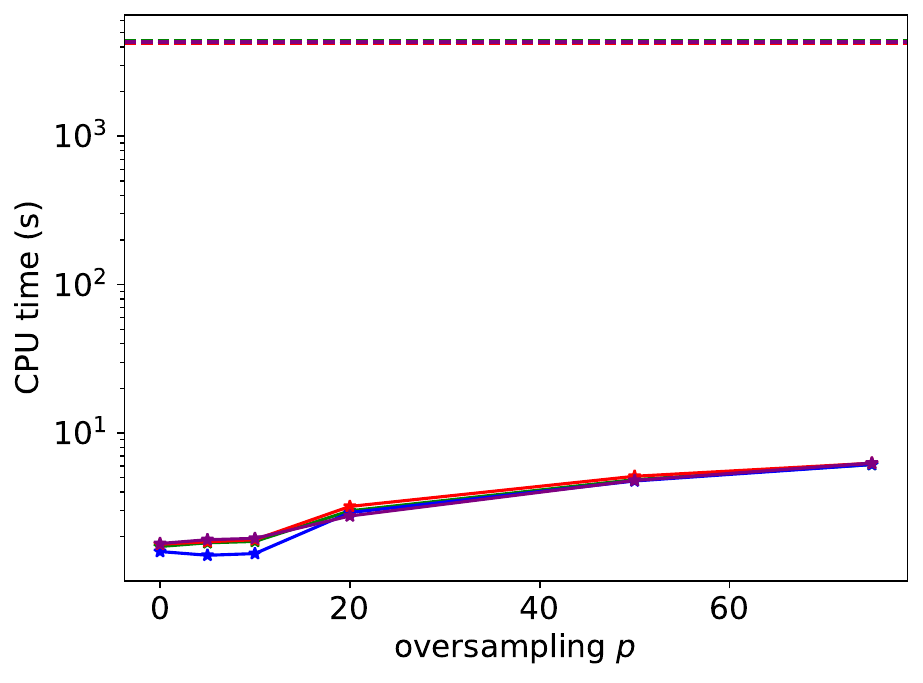}
     \label{fig:a}
 \end{subfigure}
 \hfill
 \begin{subfigure}{0.49\textwidth}
     \includegraphics[width=\textwidth]{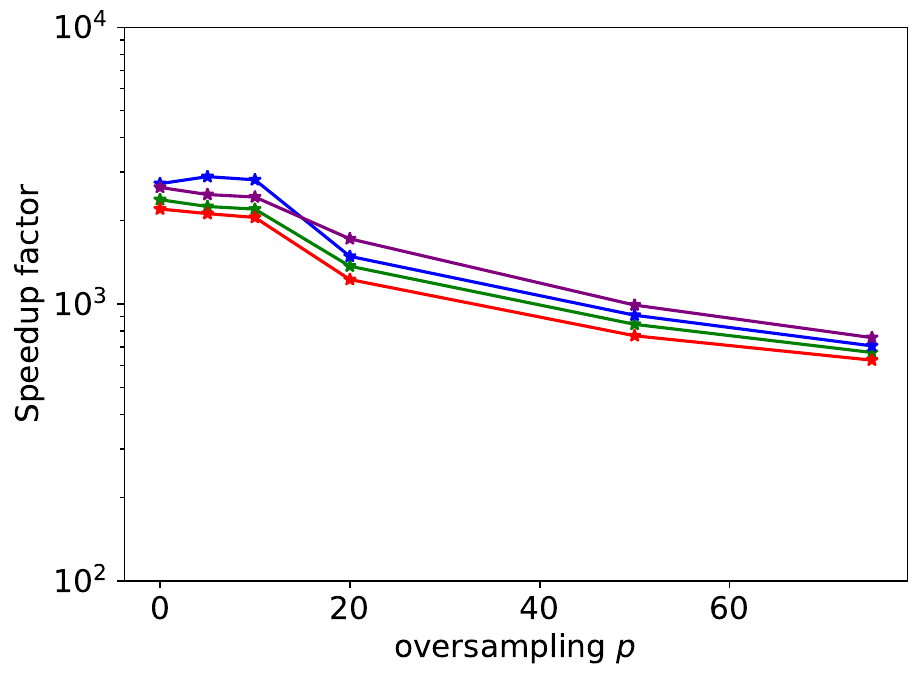}
     \label{fig:b}
 \end{subfigure}
\hfill
\begin{subfigure}{0.95\textwidth}
    \centering
    \includegraphics[width=\textwidth]{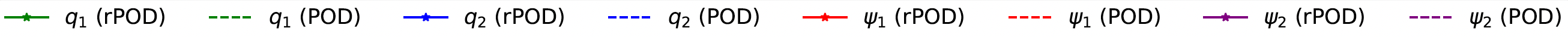}
\end{subfigure}
 \caption{Computational time needed by the POD algorithm and the rPOD algorithm (left) and the corresponding speedup factor (right) as $p$ varies.}
 \label{fig:rpod_time}
\end{figure}

The CPU time of both POD and rPOD algorithms is reported Fig.~\ref{fig:rpod_time}. The application of POD takes more than 4000 s
for each snapshot matrix, while the maximum
time required by the rPOD algorithm is roughly 6 s, 
corresponding to the snapshot matrix of $\psi_2$ with $p=75$. 
This means that rPOD 
allows for a speedup of almost 700 times. For all the tests reported in Sec.~\ref{sec:time+delta}-\ref{sec:time+delta+sigma+Fr}, we have set  
$p=75$ since it provides an excellent reconstruction of the first ten modes.

Following our previous work \cite{Besabe2024b}, we set 
$N_\Phi^r = 10$ basis functions.
As clear from Fig.~\ref{fig:singular_values}, 10 modes is nowhere near the number of modes that are needed for an accurate description of the system dynamics. In fact, for accuracy one typically retains 99.99\%
of the eigenvalue energy for each variable.
See, e.g., \cite{hesthaven2016certified}.
Tab.~\ref{tab:eig_en} reports the percentage
of eigenvalue energy retained with 10 modes for each variable as the dimension of the parameter space
increases. We see that the system dynamics is severely 
under-resolved. If we wanted to fully ($\sim$99.99\%) resolve it, that would mean retaining several hundreds of modes for each variable, which would result in computational inefficiency. Hence, we prioritize efficiency by setting $N_\Phi^r = 10$. We do expect
a limit on accuracy, especially as the parameter space
dimension increases. One way to improve accuracy is through numerical stabilization or closure models, which account for the effect of the discarded modes.
see \cite{fluids9080178} for a review ROM closures and stabilizations that are inspired by Large Eddy Simulation and \cite{Ahmed_closures2021} for a more general review of ROM closures for fluid dynamics models.

\begin{table}
    \centering
    \begin{tabular}{|c|c|c|c|c|}
         \hline
         $d+ 1$ & 2 & 3 & 4 \\
         \hline
        $q_1$ & 29\% & 29\% & 28\% \\
        $q_2$ & 26\% & 26\% & 26\% \\
        $\psi_1$ & 43\% & 39\% & 39\% \\
        $\psi_2$ & 51\% & 46\% & 46\% \\
        \hline
    \end{tabular}
    \caption{Eigenvalue energy retention when $N_\Phi^r$ = 10 modes for both $q_l$ and $\psi_l$, $l=1,2$, with increasing dimension of the physical parameter space $d$. Note that we are reporting the value of $d+1$ since time is a variable parameter in every case.}
    \label{tab:eig_en}
\end{table}

The values of the hyperparameters used for the LSTM networks  $\calmf_{q_l}$ and $\calmf_{\psi_l}$ are summarized in Table~\ref{tab:hyperparam}.
To test the performance of rPOD-LSTM ROM, we look at its reconstruction for out-of-sample set of parameters by checking the $L^2$ relative error: 
\begin{equation} \label{eq:l2-error}
    \varepsilon_{\Phi} = \frac{\norm{\widetilde{\Phi}^{\text{FOM}}-\widetilde{\Phi}^{\text{ROM}}}_{L^2(\Omega)}}{\norm{\widetilde{\Phi}^{\text{FOM}}}_{L^2(\Omega)}}.
\end{equation}
where $\tilde{\Phi}$ is the time-averaged field for $\Phi\in\{q_l,\psi_l\}$, $l =1,2$. We will consider the FOM solution as the \textit{true} solution.

\begin{table}
    \centering
    \begin{tabular}{|l|c|c|}
        \hline 
        Hyperparameters & $\calmf_{q_l}$ & $\calmf_{\psi_l}$ \\
        \hline
        Number of layers & 1 & 3 \\
        \hline
        Number of cells per layer & 100 & 50 \\
        \hline
        Batch size & 8 & 16 \\
        \hline
        Epochs & $500$ & $500$ \\
        \hline
        Activation function & $\tanh$ & $\tanh$ \\
        \hline
        Validation & $20\%$ & $20\%$ \\
        \hline
        Training : testing ratio & $1:4$ & $1:4$ \\
        \hline
        Loss function & MSE & MSE \\
        \hline
        Optimizer & Adam & Adam \\
        \hline
        Learning rate & 1E-03 & 1E-03 \\
        \hline
        Drop out probability & - & 0.1\\
        \hline
        Weight decay & 1E-05 & 1E-05 \\
        \hline
    \end{tabular}
    \caption{Hyperparameters for the LSTM network component of the rPOD-LSTM ROMs.}
    \label{tab:hyperparam}
\end{table}

\subsection{ROM for two varying parameters} \label{sec:time+delta}

\begin{figure}
    \centering
        \centering
        \begin{tabular}{cccccccccc}
            & \hspace{-0.4cm}$\delta = 0.1$ & \hspace{-0.4cm}$\delta = 0.2$ & \hspace{-0.4cm}$\delta = 0.3$ & \hspace{-0.4cm}$\delta = 0.4$ & \hspace{-0.4cm}$\delta = 0.5$ & \hspace{-0.4cm}$\delta = 0.6$ & \hspace{-0.4cm}$\delta = 0.7$ & \hspace{-0.4cm}$\delta = 0.8$ & \hspace{-0.4cm}$\delta = 0.9$ \\
            \hspace{-0.3cm}$\tildep_1$ & \hspace{-0.4cm}\includegraphics[align=c,scale = 0.6]{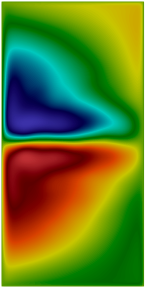} & \hspace{-0.4cm}\includegraphics[align=c,scale = 0.6]{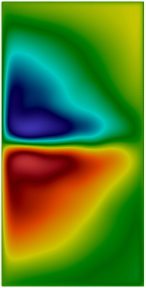} & \hspace{-0.4cm}\includegraphics[align=c,scale = 0.6]{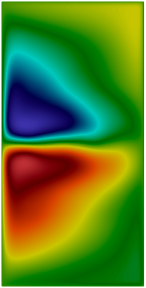} & \hspace{-0.4cm}\includegraphics[align=c,scale = 0.6]{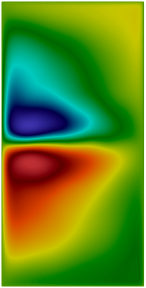} & \hspace{-0.4cm}\includegraphics[align=c,scale = 0.6]{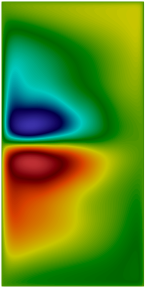} & \hspace{-0.4cm}\includegraphics[align=c,scale = 0.6]{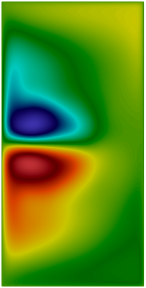} & \hspace{-0.4cm}\includegraphics[align=c,scale = 0.6]{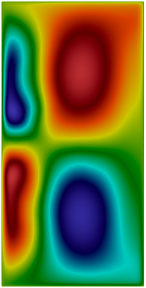} & \hspace{-0.4cm}\includegraphics[align=c,scale = 0.6]{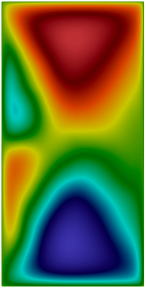} & \hspace{-0.4cm}\includegraphics[align=c,scale = 0.6]{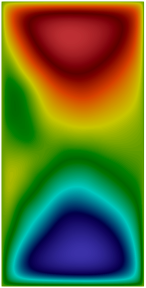} \\
            \hspace{-0.3cm}$\tildep_2$ & \hspace{-0.4cm}\includegraphics[align=c,scale = 0.6]{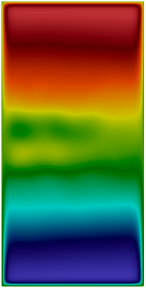} & \hspace{-0.4cm}\includegraphics[align=c,scale = 0.6]{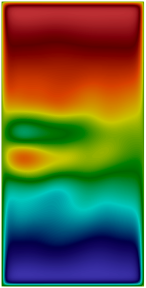} & \hspace{-0.4cm}\includegraphics[align=c,scale = 0.6]{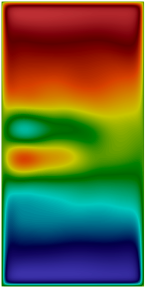} & \hspace{-0.4cm}\includegraphics[align=c,scale = 0.6]{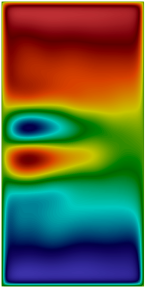} & \hspace{-0.4cm}\includegraphics[align=c,scale = 0.6]{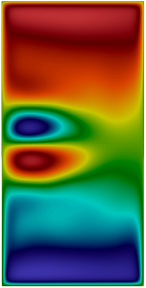} & \hspace{-0.4cm}\includegraphics[align=c,scale = 0.6]{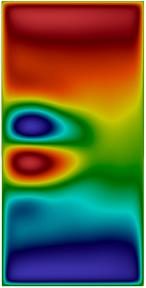} & \hspace{-0.4cm}\includegraphics[align=c,scale = 0.6]{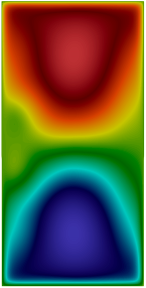} & \hspace{-0.4cm}\includegraphics[align=c,scale = 0.6]{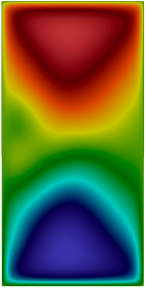} & \hspace{-0.4cm}\includegraphics[align=c,scale = 0.6]{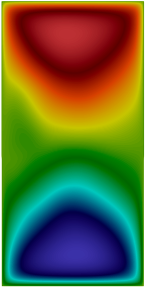}
        \end{tabular}
    \caption{Time-averaged fields $\tildep_1$ (top row) and $\tildep_2$ (bottom row) for values of $\delta$ ranging from $0.1$ (leftmost) to $0.9$ (rightmost), with $\Delta\delta = 0.1$.
    }
    \label{fig:avg-delta0.1}
\end{figure}

In this section, we will be varying two parameters: time and the aspect ratio $\delta$. All the other parameters are set as follows:
$Re = 450$, $Ro = 0.01$, $ Fr = 0.1$, and $\sigma = 0.006.$
We collect one snapshot every 0.1 time unit over the interval $[10,50]$ for sample values of
$\delta\in[0.2,0.6]$. We narrowed the interval for $\delta$ from $[0,1]$ upon a visual inspection of the time-averaged solutions. Fig.~\ref{fig:avg-delta0.1} shows the $\tildep_l$, $l = 1, 2$, for $\delta = 0.1, 0.2, \dots, 0.9$. We note that the color bar in each panel in Fig.~\ref{fig:avg-delta0.1} is set to match the respective minimum and maximum values to better show the gyre shapes and sizes.
We see that the patterns in both $\tildep_1$ and $\tildep_2$
have an abrupt change for $\delta\geq0.7$,
while for $0.1\leq\delta\leq0.6$ the shapes of the gyres remain similar, although with changes in size. For this reason, we chose to restrict the range of $\delta$ to $[0.2,0.6]$. One could 
consider the entire interval $[0,1]$ by using, e.g., local model order reduction techniques \cite{Amsallem2012, Discacciati2023, Faust2024, Vlachasa2021, Hess2019, etna_vol56_pp52-65,Hess2022}. 
However, this is beyond the scope of the current work. 

\begin{figure}[htb!]
    \centering
    \begin{subfigure}{0.45\linewidth}
        \includegraphics[width = 0.95\linewidth]{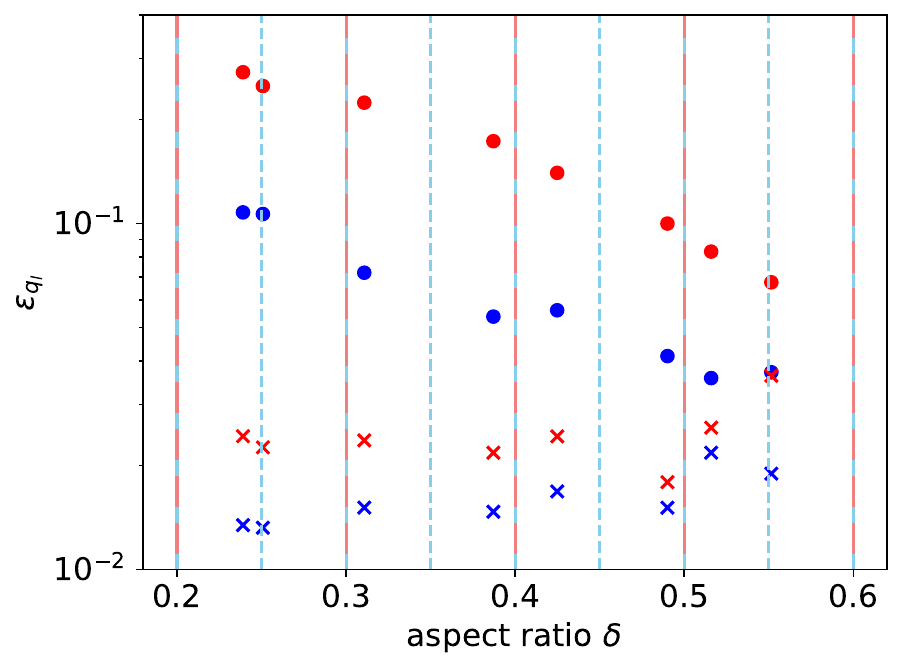}
        \caption{Potential vorticity $\tildeq_l$}
    \end{subfigure}
    \begin{subfigure}{0.45\linewidth}
        \includegraphics[width = \linewidth]{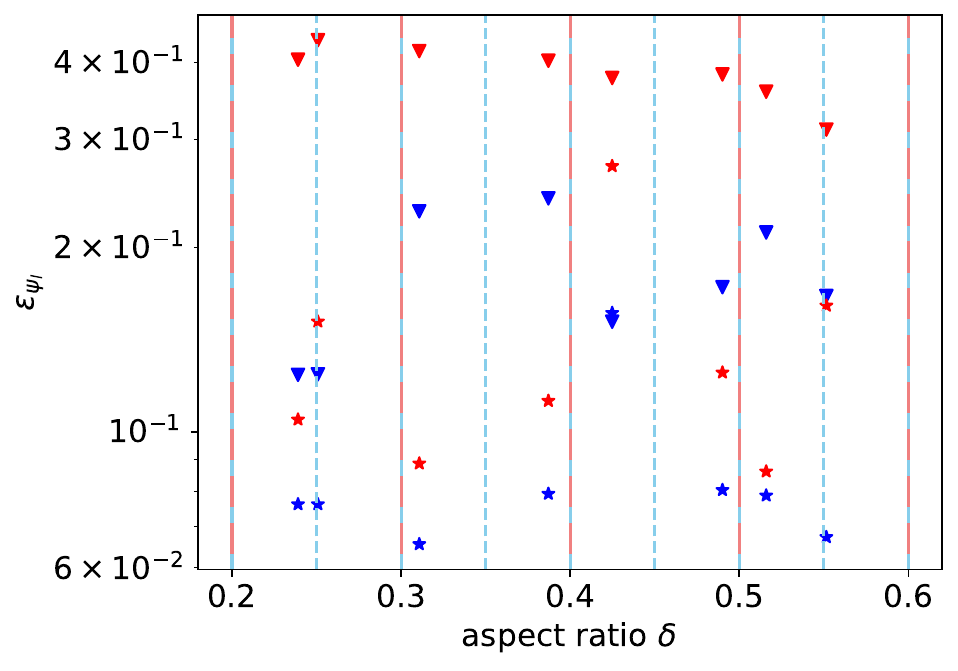}
        \caption{Stream function $\tildep_l$}
    \end{subfigure}
    \begin{subfigure}{\linewidth}
        \vspace{0.5cm}
        \includegraphics[width = \linewidth]{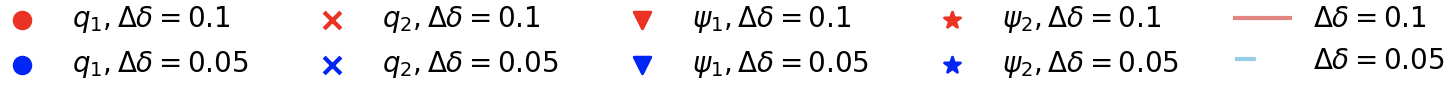}
    \end{subfigure}
    \caption{Error \eqref{eq:l2-error} for $\tildeq_l$ (left) and $\tildep_l$ (right) with uniform sampling $\Delta\delta = 0.1$ (red markers) and $\Delta\delta = 0.05$ (blue markers). Vertical lines refer to sampling values of $\delta$.}
    \label{fig:error_plot_delta}
\end{figure}

In Fig.~\ref{fig:error_plot_delta}, we show the $L^2$ relative errors $\varepsilon_\Phi$ \eqref{eq:l2-error} for $\tildeq_l$ and $\tildep_l$, $l=1,2$, using the rPOD-LSTM ROM for two different samplings for $\delta$: uniform sampling with $\Delta\delta = 0.1$ ($2005$ snapshots) and $\Delta\delta = 0.05$ ($3609$ snapshots). 
Alternatives to uniform sampling are sampling through
a greedy approach \cite{Chen2018, SIENA2023127} or at the Chebyshev nodes \cite{Pitton2017}.
From Fig.~\ref{fig:error_plot_delta},
we see that for $\Delta\delta = 0.1$ the errors are high, 
specifically around 0.2 for $\tildeq_l$ and
0.2-0.4 for $\tildep_l$. 
These errors decrease when we sample with 
$\Delta\delta = 0.05$: compare a red marker with the corresponding blue marker in Fig.~\ref{fig:error_plot_delta}.
We also note that the errors 
for $\tildep_l$ are larger for lower values of $\delta$, while the errors for $\tildeq_l$
do not increase as $\delta$ decreases. 

In the next 2 sections, we will continue sampling with $\Delta\delta = 0.05$.

\subsection{ROM for three varying parameters} \label{sec:time+delta+sigma}

In this section, we add another varying parameter: the bottom layer friction coefficient $\sigma \in [ 0.006,0.01]$. The other three parameters, i.e., $Re$, $Ro$, and $Fr$, are set like in Sec.~\ref{sec:time+delta}.

\begin{figure}[h]
    \centering
    \begin{subfigure}{0.48\linewidth}
        \includegraphics[width = \linewidth]{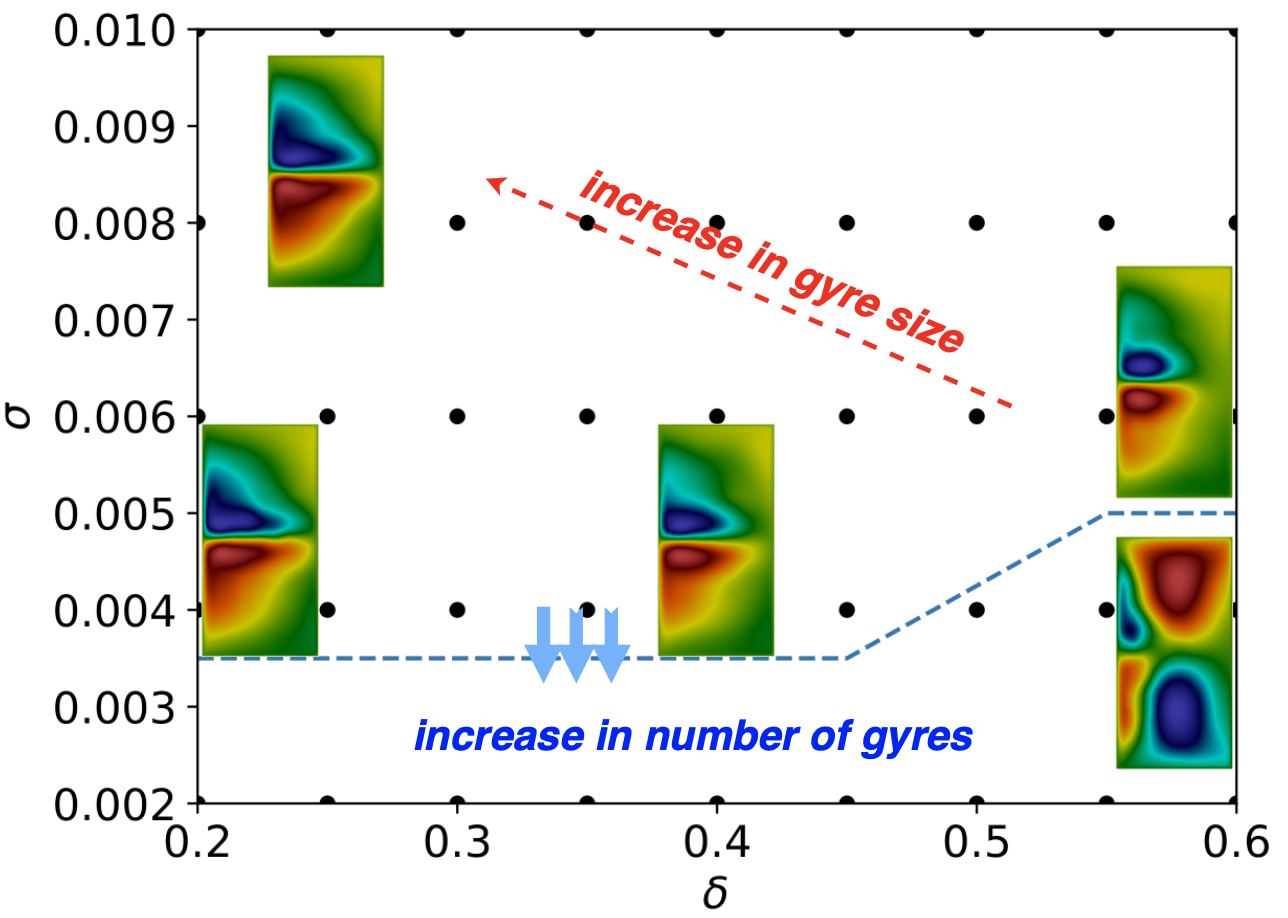}
    \end{subfigure}
    \begin{subfigure}{0.48\linewidth}
        \includegraphics[width = \linewidth]{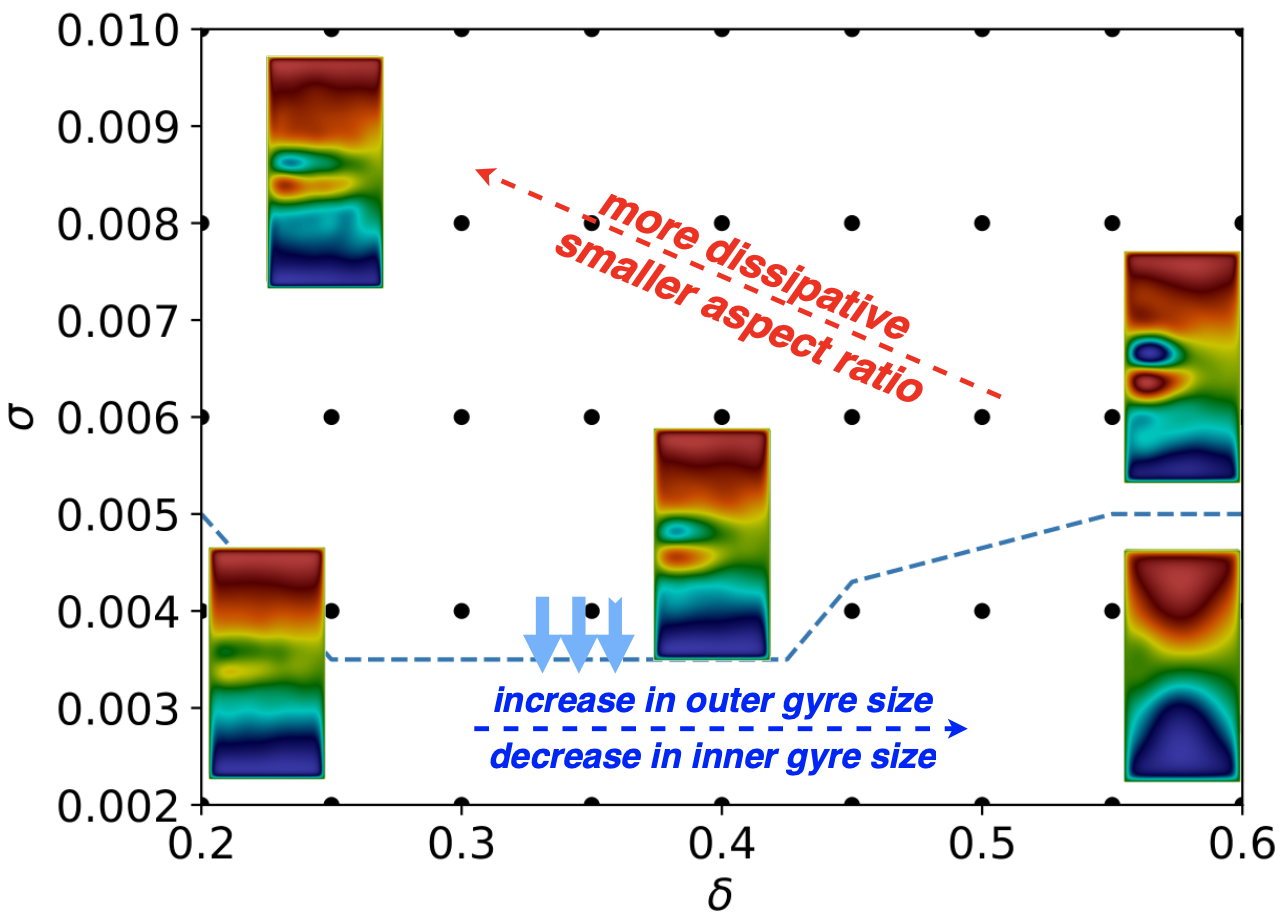}
    \end{subfigure}
    \caption{Flow pattern maps for $\tildep_1$ (left) and $\tildep_2$ (right) in the $\delta$-$\sigma$ plane.}
    \label{fig:flow_pattern_ds}
\end{figure}

Like in the case of $\delta$, we started from a larger interval for $\sigma$, i.e., $[0.002,0.01]$. Fig.~\ref{fig:flow_pattern_ds} presents a pattern maps for $\tildep_1$ (left) and $\tildep_2$ (right) in 
the $\delta$-$\sigma$ plane sampled with $\Delta\delta = 0.05$ and $\Delta \sigma = 0.001$. 
Below the blue dashed line in Fig.~\ref{fig:flow_pattern_ds} (left), we observe four gyres of equal strength in the top layer. 
Above the blue dashed line, we observe two stronger gyres at the Western boundary
that, as $\delta$ decreases and $\sigma$ increases, increase in size. 
As for the bottom layer, shown in Fig.~\ref{fig:flow_pattern_ds} (right), we see four gyres above the blue dashed line, with the the inner gyres that are long and narrow.
Below the blue dashed line, the inner gyres become faint and disappear, while the outer gyres expand. In this region, as $\delta$ increases, the outer gyres increase in size, while the inner gyres become smaller. Hence, similarly to what we have done in Sec.~\ref{sec:time+delta}, we restrict the range of $\sigma$ to $0.006\leq\sigma\leq0.010$ with uniform sampling of size $\Delta\sigma = 0.001$. 

\begin{figure}[h!]
    \centering
    \begin{subfigure}{0.45\linewidth}
        \includegraphics[width = \linewidth]{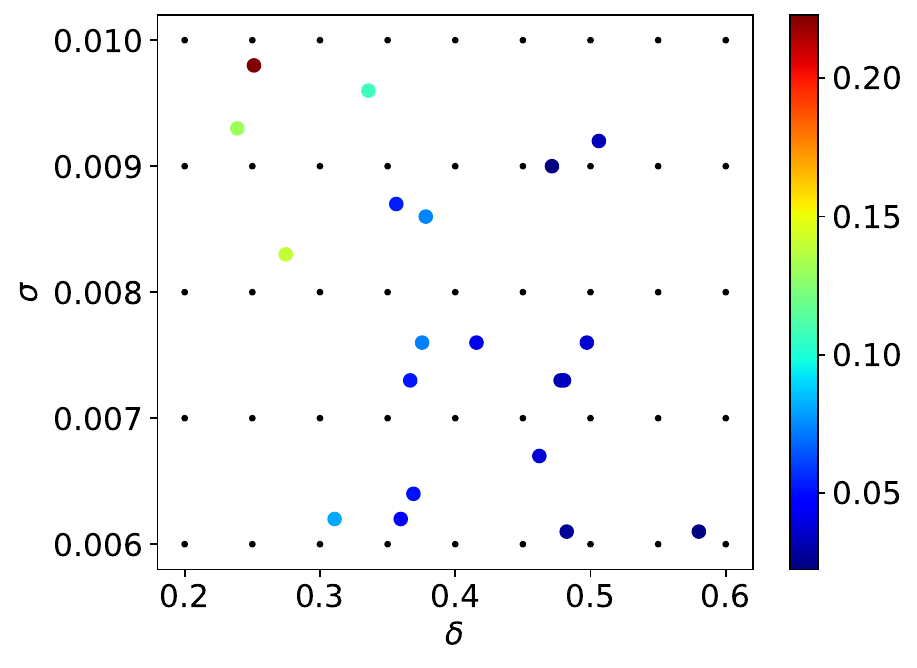}
        \caption{$\tildeq_1$}
    \end{subfigure}
    \begin{subfigure}{0.45\linewidth}
        \includegraphics[width = \linewidth]{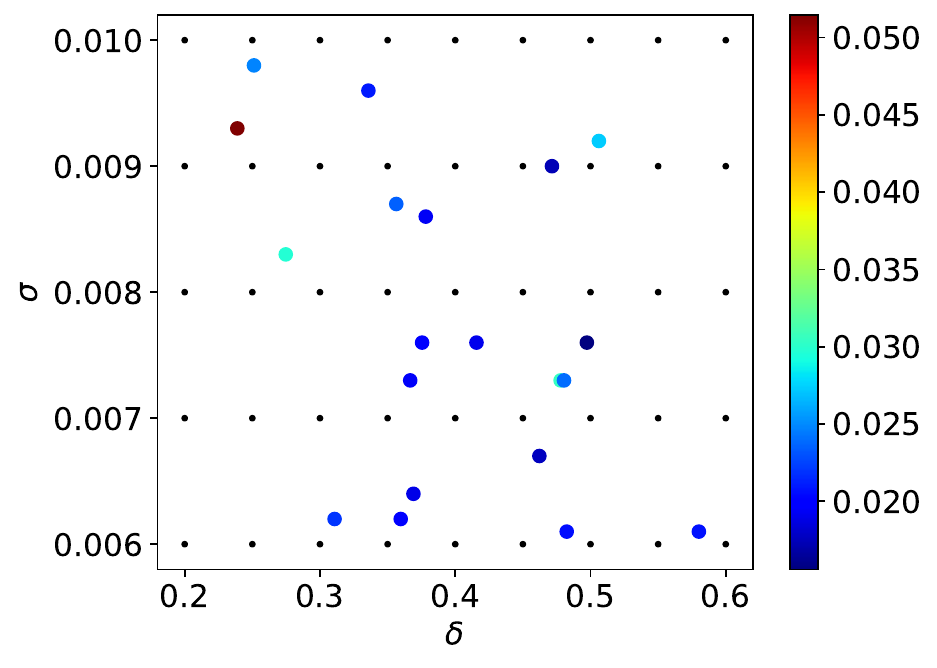}
        \caption{$\tildeq_2$}
    \end{subfigure}
    \begin{subfigure}{0.45\linewidth}
        \includegraphics[width = \linewidth]{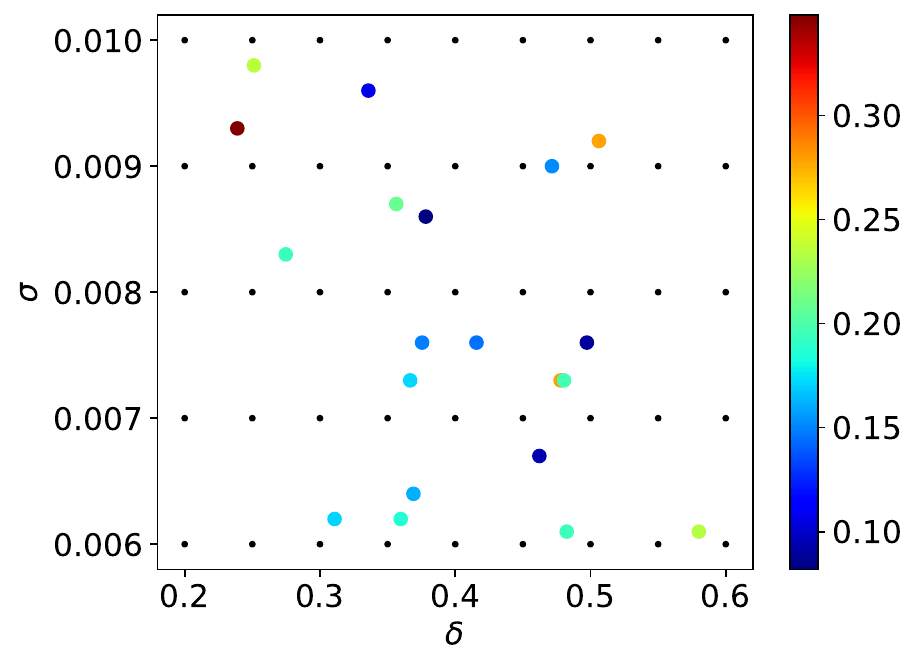}
        \caption{$\tildep_1$}
    \end{subfigure}
    \begin{subfigure}{0.45\linewidth}
        \includegraphics[width = \linewidth]{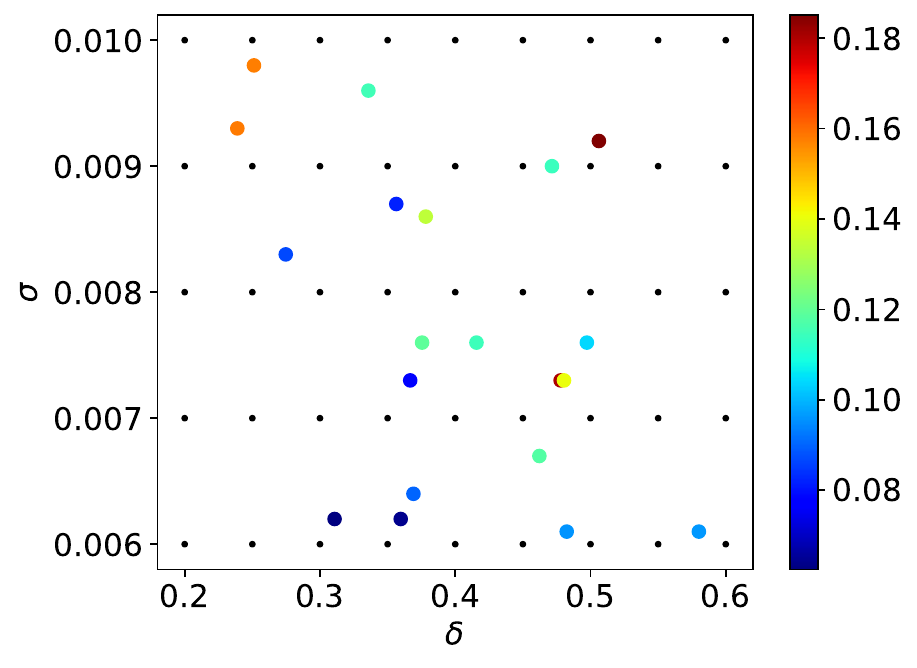}
        \caption{$\tildep_2$}
    \end{subfigure}
    \caption{Error \eqref{eq:l2-error} for the time-averaged potential vorticities $\tildeq_l$ (top row) and stream functions $\tildep_l$ (bottom row), $l = 1,2$, 
    for 20 test points. Small black dots mark the sample points.}
    \label{fig:error_plot_ds}
\end{figure}

We collect one snapshot every 0.1 time unit in interval $[10,50]$ 
for every $(\delta,\sigma)$ pair, for a total of $18,045$ 
snapshots. We apply rPOD as described in Alg.~\ref{alg:rpod} 
to generate the global basis functions. 
To assess the rPOD-LSTM ROM, we select 20 out-of-sample parameters
and compute the solution for $t \in (50,100]$, which is  
the predictive time window. The out-of-sample parameters $(\delta, \sigma)$ are chosen randomly from a uniform probability distribution over their respective ranges. 
Fig.~\ref{fig:error_plot_ds} shows error $\varepsilon_\Phi$ \eqref{eq:l2-error} 
for the time-averaged vorticity and stream function fields.
We see that errors for all the fields tend to be smaller around the center of the $0.006\leq\sigma\leq0.010$ region in the $\delta -\sigma$ plane. Notice
that the color bars are different in the different panels of Fig.~\ref{fig:error_plot_ds} because the errors for $\tildeq_l$ tend to be smaller, in average, than the errors for $\tildep_l$.
We observe larger errors (dots in dark red) for test points 
with high values of $\sigma$ and low values of $\delta$. 

\subsection{ROM for four varying parameters} \label{sec:time+delta+sigma+Fr}

In this section, we add a fourth varying parameter: the Froude number $Fr \in [0.07,0.11]$. The other two parameters, i.e., $Re$ and $Ro$, are set like in Sec.~\ref{sec:time+delta}.

\begin{figure}
    \centering
    \begin{overpic}[width=\textwidth]{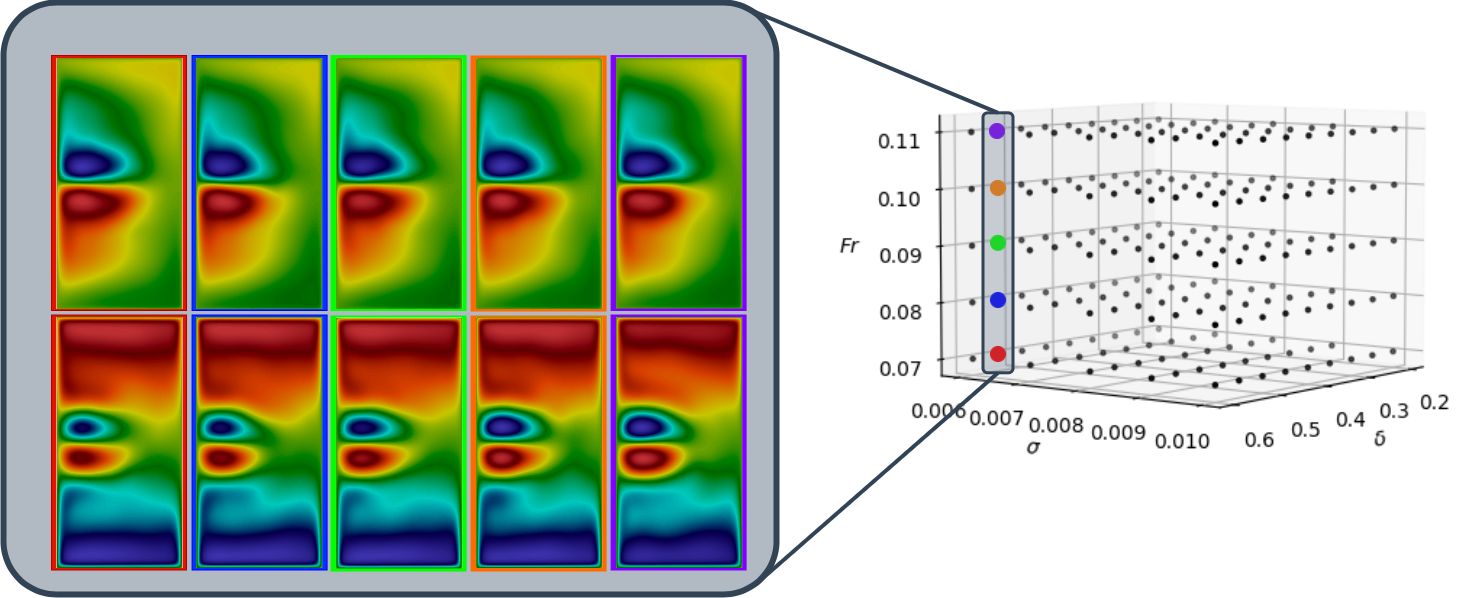}
        \put(4,38){\scriptsize{$Fr = 0.07$}}
        \put(13.5,38){\scriptsize{$Fr = 0.08$}}
        \put(23,38){\scriptsize{$Fr = 0.09$}}
        \put(32.5,38){\scriptsize{$Fr = 0.10$}}
        \put(42,38){\scriptsize{$Fr = 0.11$}}
        \put(1, 28){\scriptsize{$\tildep_1$}}
        \put(1, 10){\scriptsize{$\tildep_2$}}
    \end{overpic}
    \caption{Sample points in the $[0.2,0.6] \times [0.006,0.01] \times [0.07,0.11]$
    region in the $\delta-\sigma-Fr$ space (right) 
    and example of solutions for $\delta = 0.55$ and $\sigma = 0.006$ and varied $Fr$ (left). 
    }
    \label{fig:dsf_samples}
\end{figure}

The sampling for $\delta$ and $\sigma$ is the same used in the previous section and we consider $\Delta Fr = 0.01$. The sampling points, listed in \eqref{eq:test_delta}-\eqref{eq:test_fr}, are shown in Fig.~\ref{fig:dsf_samples} on the right. 
We visually inspected the solutions for each sampling point
and did not find drastic changes as $Fr$ is varied in $[0.07,0.11]$.
See example in Fig.~\ref{fig:dsf_samples} on the left, where we fixed $\delta = 0.55$ and $\sigma = 0.006$ and varied $Fr$. Hence, we do not further
restrict the range considered for $Fr$.
We collect one snapshot every 0.1 time unit in interval $[10,50]$ 
for every ($\delta, \sigma, Fr$) triplet, for a total of
$90,225$ snapshots.

\begin{figure}[htb!]
    \centering
    \begin{subfigure}{0.45\linewidth}
        \includegraphics[width = \linewidth]{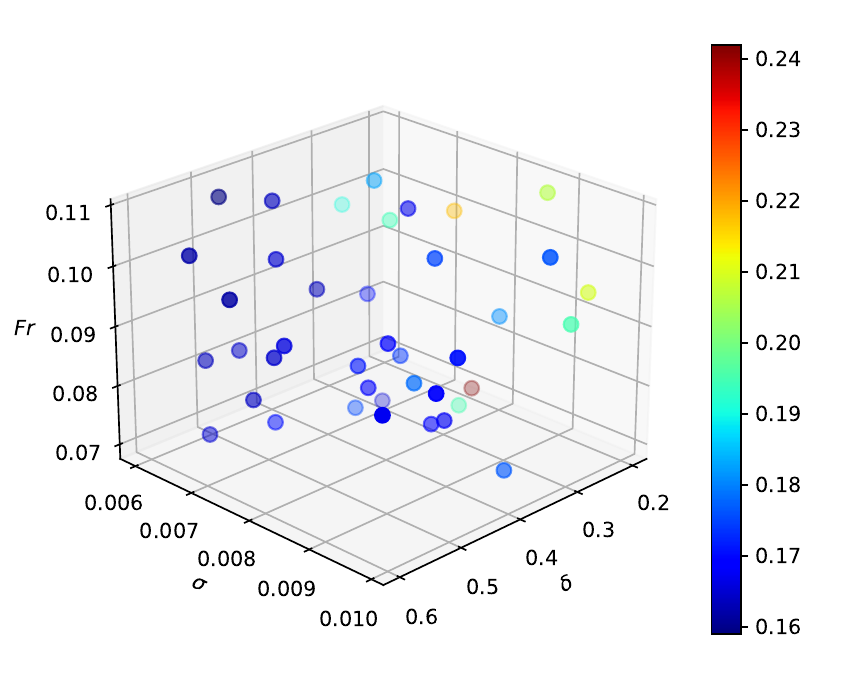}
        \caption{$\tildeq_1$}
    \end{subfigure}
    \begin{subfigure}{0.45\linewidth}
        \includegraphics[width = \linewidth]{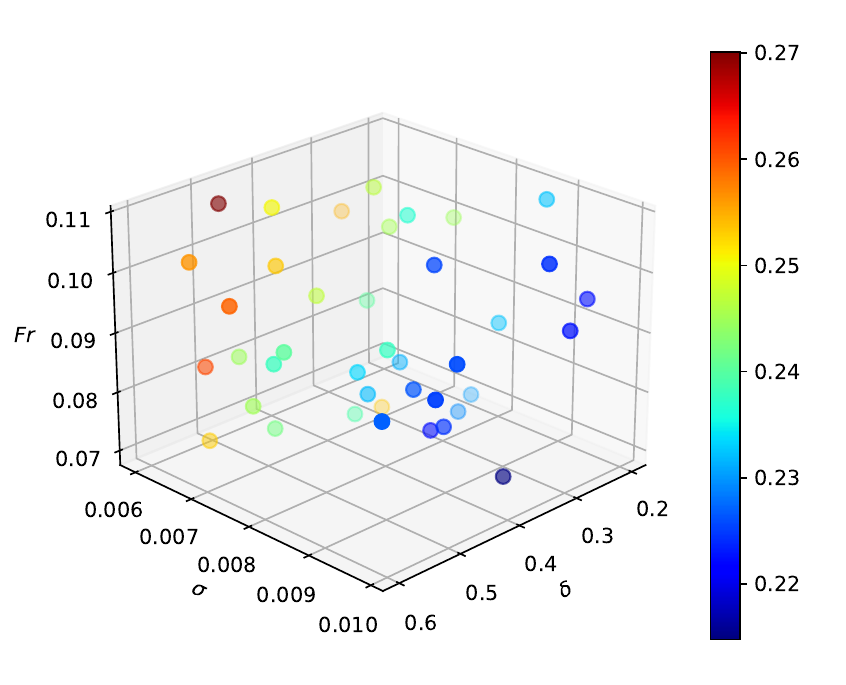}
        \caption{$\tildeq_2$}
    \end{subfigure}
    \begin{subfigure}{0.45\linewidth}
        \includegraphics[width = \linewidth]{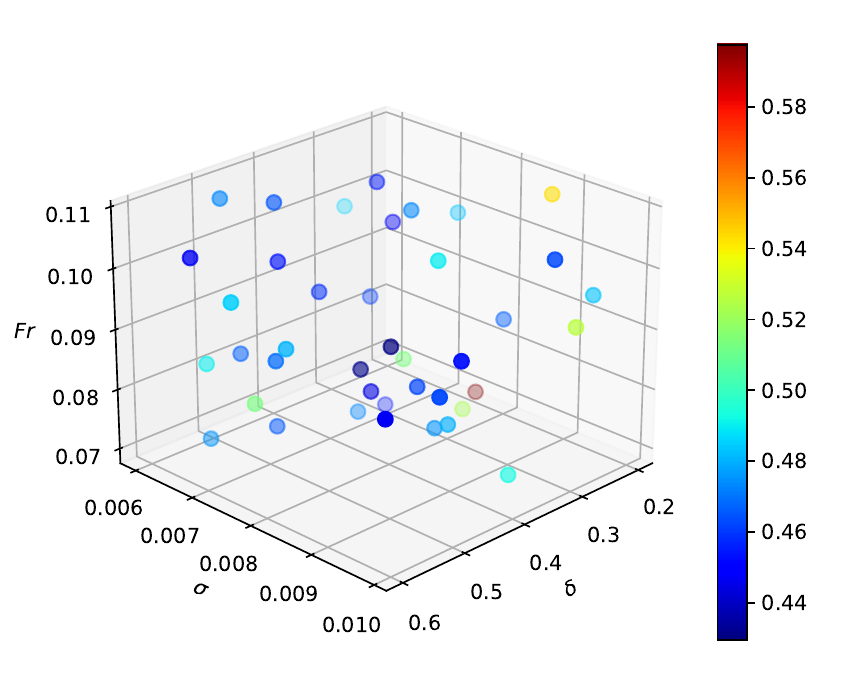}
        \caption{$\tildep_1$}
    \end{subfigure}
    \begin{subfigure}{0.45\linewidth}
        \includegraphics[width = \linewidth]{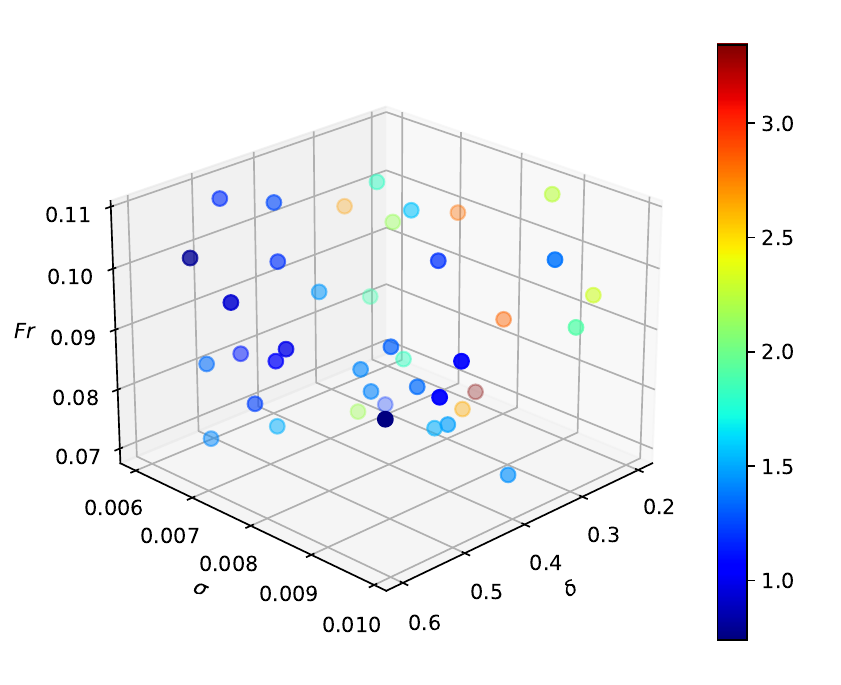}
        \caption{$\tildep_2$}
    \end{subfigure}
    \caption{
    Error \eqref{eq:l2-error} for the time-averaged potential vorticities $\tildeq_l$ (top row) and stream functions $\tildep_l$ (bottom row), $l = 1,2$, 
    for 40 test points. The sample points are shown in Fig.~\ref{fig:dsf_samples}.}
    \label{fig:error_plot_dsf}
\end{figure}

We apply the rPOD algorithm described in Alg.~\ref{alg:rpod}
to find the reduced basis functions and train the LSTM networks
as described in Sec.~\ref{sec:lstm}. 
To test the online stage of the rPOD-LSTM ROM, we generate 40 
randomly sampled parameter vectors $\bm{\mu} = (\delta, \sigma, Fr)$ 
in their respective ranges with a uniform probability distribution. 
Once again, the time interval we consider for the test is $[50,100]$, which is the predictive time window. 
Fig.~\ref{fig:error_plot_dsf} reports error \eqref{eq:l2-error} 
for the time-averaged vorticity and stream function fields.
The color bars are different in the different panels of Fig.~\ref{fig:error_plot_dsf} because the errors for $\tildeq_l$ 
tend to be smaller than the errors for $\tildep_l$.
By comparing the color bars in Fig.~\ref{fig:error_plot_ds} and \ref{fig:error_plot_dsf}, we see
that the maximum values of the errors increase 
as we go from a three-dimensional to a four-dimensional 
parameter space. This was expected because we retain 
10 modes per variable irrespective of the dimension
of the parameter space. 
In particular, from Fig.~\ref{fig:error_plot_dsf} we see that the errors for $\tildeq_2$ are larger for larger values of $\sigma$, while the errors for $\tildep_2$ are larger for smaller values of $\delta$,
which is consistent with Fig.~\ref{fig:error_plot_ds}.
A finer sampling of the intervals for $\sigma$ 
and $\delta$ could improve accuracy, as well as 
the use of a ROM closure.

\subsection{Computational time reduction} \label{sec:cpu_time}


The ROM proposed in this paper features computational savings at different levels. 

For the snapshot generation, 
we adopt a nonlinear filtering stabilization approach, denoted with 2QG-NL-$\alpha$,
that allows us to use a mesh size $h$ larger than the Munk scale.
If we were to used a ``DNS'' for the 2QGE, i.e., use a given 
eddy viscosity  and mesh size $h$ smaller than the Munk scale \eqref{eq:Munk}, 
the simulation over time interval $[0,50]$ would take around 3.7 days per sample
point in the space of physical parameters. Thanks to the 2QG-NL-$\alpha$ approach, 
that time is reduced to 3.7 hours ($\sim 24$ speedup) per sample point.

Once the snapshots are collected, we save computational time by applying
rPOD instead of POD. In the ROM with a four-dimensional parameter space, 
the application of POD takes more than 4000 s for each field of interest. 
In contrast, the application of rPOD takes only a maximum of 6 s ($\sim 700$ speedup).

During the online phase, each rPOD-LSTM ROM simulation over time interval
$[10,100]$ takes 1.5 s. Considereing that a DNS of the same time interval 
would take approximately 8.2 days, we obtain an approximate speedup of 4.7E+05.

We summarize the computational time and speedup information in Tab.~\ref{tab:cpu_time}.

\begin{table}[h]
    \centering
    \begin{tabular}{ccc}
        \hline
        Model/Technique & CPU Time & Speedup Factor\\
        \hline
        \multicolumn{3}{c}{\textit{Snapshot generation}, $t\in [10,50]$ } \\ \hline
        DNS ($h = 1/256$) & 3d 17h 48m & -  \\
        2QG-NL-$\alpha$ ($h = 1/64$) & 3h 40m & {$\sim $24} \\
        \hline
        \multicolumn{3}{c}{\textit{Reduced solution spaces construction}} \\ \hline
        POD & 4000s & -  \\
        rPOD & $\sim 6$s & {$\sim $700} \\ \hline
        \multicolumn{3}{c}{\textit{Online phase}, $t\in [10,100]$} \\ \hline
        DNS & 8d 5h & - \\
        rPOD-LSTM ROM (online) & $\sim 1.5$s & {$>$4.7E+05} \\
    \end{tabular}
    \caption{Summary of computational times and  speedup factor for three phases in our rPOD-LSTM ROM: snapshot generation, construction of the reduced basis, and online phase.}\label{tab:cpu_time}
\end{table}

\section{Conclusions} \label{sec:conclusions}

We addressed the computational challenge of simulating complex geophysical flows described by the two-layer quasi-geostrophic equations with variable parameters, by developing a data-driven ROM. Our approach combines a nonlinear filtering stabilization technique proposed in a previous work \cite{Besabe2024} to 
generate the snapshots, randomized POD (rPOD) to construct the reduced basis, and long short-term memory networks (LSTM) 
to predict the weights for the reduced basis functions
in the ROM approximation. The nonlinear filter stabilization
enables the use of a mesh size $h$ larger than the Munk scale
and thus accelerates the generation of snapshots compared to a DNS,
which required a mesh size $h$ smaller than the Munk scale. Applying rPOD to the snapshot matrices allows for an efficient computation of the reduced basis functions, achieving around 700 times speedup over deterministic POD.
The LSTM networks are trained with parameter-dependent modal coefficients of the reduced basis functions. In the online phase, the trained LSTM networks are used to autoregressively predict the modal coefficient for each field variable for out-of-sample parameters.

We assessed the performance of the rPOD-LSTM ROM with an extension of a well-known benchmark called double-gyre wind forcing experiment. We increased the dimension of the parameter space from  two (time and aspect ratio $\delta$) to three (time, $\delta$, and friction coefficient $\sigma$), and up to four (time, $\delta$, and $\sigma$, and Froude number $Fr$)
to investigate distinct dynamical regimes and their transitions with respect to parameter variations. 
To prioritize computational savings, we retain 10 modes per variable irrespective of the dimension of the parameter space. 
To assess the accuracy of the model, we selected random test points in the given parameter space, specifically 8 points 
for the two-dimensional parameter space, 20 points for the three-dimensional parameter space, and 40 points for the four-dimensional parameter space. Our results demonstrate that the model accurately predicts the time-averaged potential vorticities and stream functions of both layers (i.e., with relative $L^2$ errors ranging from 1E-02 to 1E-01) up to parameter space dimension three. The errors increase in certain regions of the parameter space when its dimension is four.
In all the cases, the computational speedup is around 4E+05 compared to a DNS for the 2QGE.

\section*{Acknowledgements} 
L. Besabe acknowledges the support of the Society for Industrial and Applied Mathematics (SIAM) for his participation in the 2024 Gene Golub SIAM Summer School, whose topic (Iterative and Randomized Methods for Large-Scale Inverse Problems) inspired this paper.
We acknowledge the support provided by PRIN “FaReX - Full and Reduced order modelling of
coupled systems: focus on non-matching methods and automatic learning” project, PNRR NGE iNEST “Interconnected Nord-Est Innovation Ecosystem” project, and INdAM-GNCS 2019–2020 projects. 

\section*{Author Declarations}
The authors have no conflicts to disclose.

\section*{Data Availability Statement}
The data that support the findings of this study are available from the corresponding author, A.~Q., upon reasonable request.

\bibliography{QGE} 

\begin{thebibliography}{10}
\expandafter\ifx\csname url\endcsname\relax
  \def\url#1{\texttt{#1}}\fi
\expandafter\ifx\csname urlprefix\endcsname\relax\def\urlprefix{URL }\fi
\expandafter\ifx\csname href\endcsname\relax
  \def\href#1#2{#2} \def\path#1{#1}\fi

\bibitem{Hoskins1975}
B.~J. Hoskins, {The geostrophic momentum approximation and the semi-geostrophic equations}, Journal of the Atmospheric Sciences 32 (1975) 233--242.

\bibitem{peter2021modelvol1}
B.~Peter, S.~Grivet-Talocia, Q.~Alfio, R.~Gianluigi, S.~Wil, L.~M. Silveira, et~al., Model Order Reduction. Volume 1: System-and Data-Driven Methods and Algorithms, De Gruyter, Berlin, Boston, 2021.
\newblock \href {https://doi.org/doi:10.1515/9783110498967} {\path{doi:doi:10.1515/9783110498967}}.

\bibitem{benner2020modelvol2}
P.~Benner, W.~Schilders, S.~Grivet-Talocia, A.~Quarteroni, G.~Rozza, L.~Miguel~Silveira, Model {{O}}rder {{R}}eduction: {{V}}olume 2: {{S}}napshot-{{B}}ased {{M}}ethods and {{A}}lgorithms, De Gruyter, 2020.

\bibitem{benner2020modelvol3}
P.~Benner, W.~Schilders, S.~Grivet-Talocia, A.~Quarteroni, G.~Rozza, L.~Miguel~Silveira, Model order reduction: volume 3 applications, De Gruyter, 2020.

\bibitem{hesthaven2016certified}
J.~S. Hesthaven, G.~Rozza, B.~Stamm, Certified reduced basis methods for parametrized partial differential equations, Vol. 590, Springer, 2016.

\bibitem{malik2017reduced}
M.~H. Malik, Reduced {{O}}rder {{M}}odeling for smart grids' simulation and optimization, Ph.D. thesis, {\'E}cole centrale de Nantes; Universitat polit{\'e}cnica de Catalunya (2017).

\bibitem{rozza2008reduced}
G.~Rozza, D.~B.~P. Huynh, A.~T. Patera, Reduced basis approximation and a posteriori error estimation for affinely parametrized elliptic coercive partial differential equations: application to transport and continuum mechanics, Archives of Computational Methods in Engineering 15~(3) (2008) 229--275.

\bibitem{Berkooz1993}
G.~Berkooz, P.~Holmes, J.~L. Lumley, {The proper orthogonal decomposition in the analysis of turbulent flows}, Annual Review of Fluid Mechanics 25 (1993) 539--575.

\bibitem{Schmid2010}
P.~Schmid, {Dynamic mode decomposition of numerical and experimental data}, Journal of Fluid Mechanics 656 (2010) 5--28.

\bibitem{Casenave2019}
F.~Casenave, N.~Akkari, F.~Bordeu, C.~Rey, D.~Ryckelynck, {A nonintrusive distributed reduced-order modeling framework for nonlinear structural mechanics—Application to elastoviscoplastic computations}, International Journal for Numerical Methods in Engineering 121 (2019) 32--53.

\bibitem{Guo2018}
M.~Guo, J.~S. Hesthaven, {Reduced order modeling for nonlinear structural analysis using Gaussian process regression}, Computer Methods in Applied Mechanics and Engineering 341 (2018) 807--826.

\bibitem{Lin1995}
R.~Lin, M.~Lim, {Structural sensitivity analysis via reduced-order analytical model}, Computer Methods in Applied Mechanics and Engineering 121 (1995) 345--359.

\bibitem{Oliver2017}
J.~Oliver, M.~Caicedo, A.~Huespe, J.~Hernández, E.~Roubin, {Reduced order modeling strategies for computational multiscale fracture}, Computer Methods in Applied Mechanics and Engineering 313 (2017) 560--595.

\bibitem{Lucas2019}
L.~O. M\"{u}ller, A.~Caiazzo, P.~J. Blanco, {Reduced-Order Unscented Kalman Filter With Observations in the Frequency Domain: Application to Computational Hemodynamics}, IEEE Transactions on Biomedical Engineering 66 (2019) 1269--1276.

\bibitem{Pfaller2020}
M.~R. Pfaller, M.~C. Varona, C.~B. J.~Lang, W.~A. Wall, {Using parametric model order reduction for inverse analysis of large nonlinear cardiac simulations}, International Journal for Numerical Methods in Biomedical Engineering 36 (2020) e3320.

\bibitem{Pfaller2022}
M.~R. Pfaller, J.~Pham, A.~Verma, L.~Pegolotti, N.~M. Wilson, D.~W. Parker, W.~Yang, A.~L. Marsden, {Automated generation of 0D and 1D reduced-order modelsof patient-specific blood flow}, International Journal for Numerical Methods in Biomedical Engineering 38 (2022) e3639.

\bibitem{Itu2012}
L.~Itu, P.~Sharma, V.~Mihalef, A.~Kamen, C.~Suciu, D.~Lomaniciu, {A patient-specific reduced-order model for coronary circulation}, in: 2012 9th IEEE International Symposium on Biomedical Imaging (ISBI), 2012.

\bibitem{Zhang2020}
X.~Zhang, D.~Wu, F.~Miao, H.~Liu, Y.~Li, {Personalized Hemodynamic Modeling of the Human Cardiovascular System: A Reduced-Order Computing Model}, IEEE Transactions on Biomedical Engineering 67 (2020) 2754--2764.

\bibitem{Bryan1963}
K.~Bryan, A numerical investigation of a nonlinear model of a wind-driven ocean, Journal of Atmospheric Sciences 20~(6) (1963) 594 -- 606.
\newblock \href {https://doi.org/10.1175/1520-0469(1963)020<0594:ANIOAN>2.0.CO;2} {\path{doi:10.1175/1520-0469(1963)020<0594:ANIOAN>2.0.CO;2}}.

\bibitem{Gates1968}
W.~L. Gates, {A Numerical Study of Transient Rossby Waves in a Wind-Driven Homogeneous Ocean}, Journal of Atmospheric Sciences 25~(1) (1968) 3 -- 22.
\newblock \href {https://doi.org/10.1175/1520-0469(1968)025<0003:ANSOTR>2.0.CO;2} {\path{doi:10.1175/1520-0469(1968)025<0003:ANSOTR>2.0.CO;2}}.

\bibitem{Holland1975}
W.~R. Holland, L.~B. Lin, On the generation of mesoscale eddies and their contribution to the oceanicgeneral circulation. ii. a parameter study, Journal of Physical Oceanography 5~(4) (1975) 658 -- 669.
\newblock \href {https://doi.org/10.1175/1520-0485(1975)005<0658:OTGOME>2.0.CO;2} {\path{doi:10.1175/1520-0485(1975)005<0658:OTGOME>2.0.CO;2}}.

\bibitem{Berloff1999}
P.~S. Berloff, J.~C. McWilliams, Large-scale, low-frequency variability in wind-driven ocean gyres, Journal of Physical Oceanography 29~(8) (1999) 1925 -- 1949.
\newblock \href {https://doi.org/10.1175/1520-0485(1999)029<1925:LSLFVI>2.0.CO;2} {\path{doi:10.1175/1520-0485(1999)029<1925:LSLFVI>2.0.CO;2}}.

\bibitem{BERLOFF_KAMENKOVICH_PEDLOSKY_2009}
P.~Berloff, I.~Kamenovich, J.~Pedlosky, A mechanism of formation of multiple zonal jets in the oceans, Journal of Fluid Mechanics 628 (2009) 395–425.
\newblock \href {https://doi.org/10.1017/S0022112009006375} {\path{doi:10.1017/S0022112009006375}}.

\bibitem{Tanaka2010}
Y.~Tanaka, K.~Akitomo, Alternating zonal flows in a two-layer wind-driven ocean, J. Oceanogr. 66 (2010) 475–487.

\bibitem{Besabe2024}
L.~Besabe, M.~Girfoglio, A.~Quaini, G.~Rozza, {Linear and nonlinear filtering for a two-layer quasi-geostrophic ocean model}, Applied Mathematics and Computation 488 (2025) 129121.

\bibitem{Erichson2019}
N.~B. Erichson, L.~Mathelin, J.~N. Kutz, S.~L. Brunton, {Randomized Dynamic Mode Decomposition}, SIAM Journal on Applied Dynamical Systems 18 (2019) 1867 -- 1891.

\bibitem{Bistrian2018}
D.~A. Bistrian, I.~M. Navon, {Efficiency of randomised dynamic mode decomposition for reduced order modelling}, International Journal of Computational Fluid Dynamics 32 (2018) 88--103.

\bibitem{Tissot2014}
G.~Tissot, L.~Cordier, N.~Benard, B.~Noack, {Model reduction using Dynamic Mode Decomposition}, Comptes Rendus Mecanique 342 (2014) 410–416.

\bibitem{Buhr2018}
A.~Buhr, K.~Smetana, {Randomized Local Model Order Reduction}, SIAM Journal on Scientific Computing 40 (2018) A2120--A2151.

\bibitem{Saibaba2020}
A.~K. Saibaba, {Randomized Discrete Empirical Interpolation Method for Nonlinear Model Reduction}, SIAM Journal on Scientific Computing 42 (2020) A1582--A1608.

\bibitem{Peherstorfer2014}
B.~Peherstorfer, D.~Butnaru, K.~Willcox, H.-J. Bungartz, {Localized Discrete Empirical Interpolation Method}, SIAM Journal on Scientific Computing 36 (2014) A168--A192.

\bibitem{Yu2020}
D.~Yu, S.~Chakravorty, {A randomized balanced proper orthogonal decomposition technique}, Journal of Computational and Applied Mathematics 368 (2020) 112540.

\bibitem{Rowley2005}
C.~Rowley, {Model Reduction for Fluids, using balanced proper orthogonal decomposition}, International Journal of Bifurcation and Chaos 15 (2005) 997--1013.

\bibitem{Alla2019}
A.~Alla, J.~N. Kutz, {Randomized model order reduction}, Advances in Computational Mathematics 45 (2019) 1251--1271.

\bibitem{Rajaram2020}
D.~Rajaram, C.~Perron, T.~G. Puranik, D.~N. Mavris, {Randomized Algorithms for Non-Intrusive Parametric Reduced Order Modeling}, AIAA Journal 58 (2020) 5389--5407.

\bibitem{Bach2019}
C.~Bach, D.~Ceglia, L.~Song, F.~Duddeck, {Randomized low-rank approximation methods for projection-based model order reduction of large nonlinear dynamical problems}, International Journal for Numerical Methods in Engineering 118 (2019) 209--241.

\bibitem{Besabe2024b}
L.~Besabe, M.~Girfoglio, A.~Quaini, G.~Rozza, {Data-driven reduced order modeling of a two-layer quasi-geostrophic ocean model}, Results in Engineering 25 (2025) 103691.

\bibitem{Rahman-2019}
S.~M. Rahman, S.~Pawar, O.~San, T.~Iliescu, Nonintrusive reduced order modeling framework for quasigeostrophic turbulence, Physical Review E 100 (2019) 053306.

\bibitem{Girfoglio2023}
M.~Girfoglio, A.~Quaini, G.~Rozza, A linear filter regularization for {POD-based} reduced-order models of the quasi-geostrophic equations, Comptes Rendus. M\'ecanique 351 (2023) 1--21.

\bibitem{fluids9080178}
A.~Quaini, O.~San, A.~Veneziani, T.~Iliescu, \href{https://www.mdpi.com/2311-5521/9/8/178}{Bridging large eddy simulation and reduced-order modeling of convection-dominated flows through spatial filtering: Review and perspectives}, Fluids 9~(8) (2024).
\newblock \href {https://doi.org/10.3390/fluids9080178} {\path{doi:10.3390/fluids9080178}}.
\newline\urlprefix\url{https://www.mdpi.com/2311-5521/9/8/178}

\bibitem{Marshall1997}
J.~Marshall, C.~Hill, L.~Perelman, A.~Adcroft, Hydrostatic, quasi-hydrostatic, and nonhydrostatic ocean modeling, Journal of Geophysical Research: Oceans 102~(C3) (1997) 5733--5752.
\newblock \href {https://doi.org/https://doi.org/10.1029/96JC02776} {\path{doi:https://doi.org/10.1029/96JC02776}}.

\bibitem{Chassignet1998}
T.~M. Ozgokmen, E.~P. Chassignet, Emergence of inertial gyres in a two-layer quasigeostrophic ocean model, Journal of Physical Oceanography 28~(3) (1998) 461 -- 484.
\newblock \href {https://doi.org/10.1175/1520-0485(1998)028<0461:EOIGIA>2.0.CO;2} {\path{doi:10.1175/1520-0485(1998)028<0461:EOIGIA>2.0.CO;2}}.

\bibitem{DiBattista2001}
M.~DiBattista, A.~Majda, Equilibrium statistical predictions for baroclinic vortices: The role of angular momentum, Theoret. Comput. Fluid Dynamics 14 (2001) 293–322.

\bibitem{San2012}
O.~San, A.~Staples, T.~Iliescu, Approximate deconvolution large eddy simulation of a stratified two-layer quasigeostrophic ocean model, Ocean Modelling 63 (2012) 1--20.

\bibitem{Salmon1978}
R.~Salmon, Two-layer quasi-geostrophic turbulence in a simple special case, Geophysical \& Astrophysical Fluid Dynamics 10 (1978) 25--52.

\bibitem{Medjo2000}
T.~T. Medjo, Numerical simulations of a two-layer quasi-geostrophic equation of the ocean, SIAM Journal of Numerical Analysis 37 (2000) 2005--2022.

\bibitem{Fandry1984}
C.~B. Fandry, L.~M. Leslie, A two-layer quasi-geostrophic model of summer trough formation in the australian subtropical easterlies, Journal of the Atmospheric Sciences 41 (1984) 807--818.

\bibitem{Mu1994}
M.~Mu, Z.~Qingcun, T.~G. Shepherd, L.~Yongming, Nonlinear stability of multilayer quasi-geostrophic flow, Journal of Fluid Dynamics 264 (1994) 165--184.

\bibitem{Ikeda1981}
M.~Ikeda, Meanders and detached eddies of a strong eastward-flowing jet using a two-layer quasi-geostrophic model, Journal of Physical Oceanography 11 (1981) 526–540.

\bibitem{Zalesny2022}
V.~B. Zalesny, {Variational method for solving the Qqasi-geostrophic circulation problem in a two-Layer ocean}, Izvestiya, Atmospheric and Oceanic Physics 58 (2022) 423–432.

\bibitem{Cummins1992}
P.~Cummins, {Inertial gyres in decaying and forced geostrophic turbulence}, Journal of Marine Research 50 (1992) 545--566.

\bibitem{Ozgok1998}
T.~M. \"{O}zg\"{o}kmen, E.~P. Chassignet, {Emergence of Inertial Gyres in a Two-Layer Quasigeostrophic Ocean Model}, Journal of Physical Oceanography 28 (1998) 461–484.

\bibitem{Girfoglio_JCAM2023}
M.~Girfoglio, A.~Quaini, G.~Rozza, {A novel Large Eddy Simulation model for the Quasi-Geostrophic equations in a Finite Volume setting}, Journal of Computational and Applied Mathematics 418 (2023) 114656.

\bibitem{GEA}
{GEA - Geophysical and Environmental Applications}, \url{https://github.com/GEA-Geophysical-and-Environmental-Apps/GEA} (2023).

\bibitem{GirfoglioFVCA10}
M.~Girfoglio, A.~Quaini, G.~Rozza, {GEA: A New Finite Volume-Based Open Source Code for the Numerical Simulation of Atmospheric and Ocean Flows}, in: Finite Volumes for Complex Applications X—Volume 2, Hyperbolic and Related Problems, 2023, pp. 151--159.

\bibitem{Weller1998}
H.~G. Weller, G.~Tabor, H.~Jasak, C.~Fureby, A tensorial approach to computational continuum mechanics using object-oriented techniques, Computers in Physics 12 (1998) 620--631.

\bibitem{Martinsson2016}
P.-G. Martinsson, {Randomized methods for matrix computations}, https://arxiv.org/pdf/1607.01649 (2016).

\bibitem{Gu2015}
M.~Gu, {Subspace Iteration Randomization and Singular Value Problems}, SIAM Journal on Scientific Computing 37 (2015) A1139--A1173.

\bibitem{HMT2011}
N.~Halko, P.~Martinsson, J.~Tropp, {Finding Structure with Randomness: Probabilistic Algorithms forConstructing ApproximateMatrix Decompositions}, SIAM Review 53 (2011) 217--288.

\bibitem{Golzar2023}
M.~Golzar, M.~K. Moayyedi, F.~Fotouhi, A surrogate non-intrusive reduced order model of quasi-geostrophic turbulence dynamics based on a combination of lstm and different approaches of dmd, Journal of Turbulence 24 (2023) 474--505.

\bibitem{LSTM1997}
S.~Hochreiter, J.~Schmidhuber, {Long Short-Term Memory}, Neural Computation 9 (1997) 1735–1780.

\bibitem{Nadiga2001}
B.~Nadiga, L.~Margolin, Dispersive-dissipative eddy parameterization in a barotropic model, Journal of Physical Oceanography 31 (2001) 2525--2531.

\bibitem{Holm2003}
D.~Holm, B.~Nadiga, Modeling mesoscale turbulence in the barotropic double-gyre circulation, Journal of Physical Oceanography 33 (2003) 2355--2365.

\bibitem{Greathbatch2000}
R.~Greatbatch, B.~Nadiga, Four-gyre circulation in a barotropic model with double-gyre wind forcing, Journal of Physical Oceanography 30 (2000) 1461--1471.

\bibitem{Monteiro2014}
I.~Monteiro, C.~Carolina, Improving numerical accuracy in a regularized barotropic vorticity model of geophysical flow, International Journal of Numerical Analysis and Modelling, Series B 5 (2014) 317--338.

\bibitem{Monteiro2015}
I.~Monteiro, C.~Manica, L.~Rebholz, Numerical study of a regularized barotropic vorticity model of geophysical flow, Numerical Methods for Partial Differential Equations 31 (2015) 1492--1514.

\bibitem{QGE-review}
Z.~Mou, C.and~Wang, D.~Wells, X.~Xie, T.~Iliescu, Reduced order models for the quasi-geostrophic equations: A brief survey, Fluids 6 (2020) 16.
\newblock \href {https://doi.org/10.3390/fluids6010016} {\path{doi:10.3390/fluids6010016}}.

\bibitem{San2015}
O.~San, T.~Iliescu, A stabilized proper orthogonal decomposition reduced-order model for large scale quasigeostrophic ocean circulation, Advances in Computational Mathematics 41 (2015) 1289--1319.

\bibitem{Ahmed_closures2021}
S.~E. Ahmed, S.~Pawar, O.~San, A.~Rasheed, T.~Iliescu, B.~R. Noack, \href{https://doi.org/10.1063/5.0061577}{{On closures for reduced order models - A spectrum of first-principle to machine-learned avenues}}, Physics of Fluids 33~(9) (2021) 091301.
\newblock \href {http://arxiv.org/abs/https://pubs.aip.org/aip/pof/article-pdf/doi/10.1063/5.0061577/13869240/091301\_1\_online.pdf} {\path{arXiv:https://pubs.aip.org/aip/pof/article-pdf/doi/10.1063/5.0061577/13869240/091301\_1\_online.pdf}}, \href {https://doi.org/10.1063/5.0061577} {\path{doi:10.1063/5.0061577}}.
\newline\urlprefix\url{https://doi.org/10.1063/5.0061577}

\bibitem{Amsallem2012}
D.~Amsallem, M.~J. Zahr, C.~Farhat, {Nonlinear model order reduction based on local reduced-order bases}, International Journal for Numerical Methods in Engineering 92 (2012) 835--916.

\bibitem{Discacciati2023}
N.~Discacciati, J.~S. Hesthaven, {Localized model order reduction and domain decomposition methods for coupled heterogeneous systems}, International Journal for Numerical Methods in Engineering 124 (2023) 3964 -- 3996.

\bibitem{Faust2024}
E.~Faust, L.~Scheunemann, {Nonlinear model order reduction using local basis methods on RVEs: A parameter study and comparison of different variants}, in: 94th Annual Meeting of the International Association of Applied Mathematics and Mechanics (GAMM), 2024.

\bibitem{Vlachasa2021}
K.~Vlachasa, K.~Tatsisa, K.~Agathosb, A.~R. Brinkc, E.~Chatzi, {A local basis approximation approach for nonlinear parametric model order reduction}, Journal of Sound and Vibration 502 (2021) 116055.

\bibitem{Hess2019}
M.~Hess, A.~Alla, A.~Quaini, G.~Rozza, M.~Gunzburger, {A localized reduced-order modeling approach for PDEs with bifurcating solutions}, Computer Methods in Applied Mechanics and Engineering 351 (2019) 379--403.

\bibitem{etna_vol56_pp52-65}
M.~W. Hess, A.~Quaini, G.~Rozza, {A comparison of reduced-order modeling approaches using artificial neural networks for PDEs with bifurcating solutions}, Electron. Trans. Numer. Anal. 56 (2022) 52--65.
\newblock \href {https://doi.org/10.1553/etna_vol56s52} {\path{doi:10.1553/etna_vol56s52}}.

\bibitem{Hess2022}
M.~W. Hess, A.~Quaini, G.~Rozza, {Data-Driven Enhanced Model Reduction for Bifurcating Models in Computational Fluid Dynamics}, in: eccomas2022, 2022.
\newblock \href {https://doi.org/https://www.scipedia.com/public/Rozza_2022a} {\path{doi:https://www.scipedia.com/public/Rozza_2022a}}.

\bibitem{Chen2018}
W.~Chen, J.~S. Hesthaven, B.~Junqiang, Y.~Qiu, Z.~Yang, Y.~Tihao, {Greedy Nonintrusive Reduced Order Model for Fluid Dynamics}, AIAA Journal 56 (2018) 4927 -- 4943.

\bibitem{SIENA2023127}
P.~Siena, M.~Girfoglio, G.~Rozza, \href{https://www.sciencedirect.com/science/article/pii/B9780323899673000081}{Chapter 6 - an introduction to pod-greedy-galerkin reduced basis method}, in: F.~Chinesta, E.~Cueto, Y.~Payan, J.~Ohayon (Eds.), Reduced Order Models for the Biomechanics of Living Organs, Biomechanics of Living Organs, Academic Press, 2023, pp. 127--145.
\newblock \href {https://doi.org/https://doi.org/10.1016/B978-0-32-389967-3.00008-1} {\path{doi:https://doi.org/10.1016/B978-0-32-389967-3.00008-1}}.
\newline\urlprefix\url{https://www.sciencedirect.com/science/article/pii/B9780323899673000081}

\bibitem{Pitton2017}
G.~Pitton, A.~Quaini, G.~Rozza, {Computational reduction strategies for the detection of steady bifurcations in incompressible fluid-dynamics: Applications to Coanda effect in cardiology}, Journal of Computational Physics 344 (2017) 534--557.

\end{thebibliography}

\end{document}